\documentclass[final,3p,times]{elsarticle}

\usepackage{amsmath,amssymb,amsfonts,amsthm}
\usepackage{bm}
\usepackage{mathtools}
\usepackage{mathrsfs}
\usepackage{algorithm}
\usepackage{algpseudocode}
\usepackage{booktabs}
\usepackage{graphicx}
\usepackage{float}
\usepackage{placeins}
\usepackage{microtype}
\usepackage{tikz}
\usetikzlibrary{arrows.meta,positioning}
\usepackage[colorlinks=true,linkcolor=blue,citecolor=blue,urlcolor=blue]{hyperref}
\usepackage{comment}

\usepackage{graphicx}
\usepackage{xcolor}
\usepackage{subcaption}
\usepackage[textsize=footnotesize]{todonotes}


\newtheorem{proposition}{Proposition}

\newcommand{\R}{\mathbb{R}}
\newcommand{\Ncal}{\mathcal{N}}
\newcommand{\Acal}{\mathcal{A}}
\newcommand{\Bcal}{\mathcal{B}}
\newcommand{\Ccal}{\mathcal{C}}

\newcommand{\Oh}{\widehat{\Omega}}

\newcommand{\Phihat}{\widehat{\Phi}_K}
\newcommand{\Phat}{\widehat{P}_K}
\newcommand{\x}{\bm{x}}
\newcommand{\xiq}{\bm{\xi}}
\newcommand{\cp}{\bm{c}_p}
\newcommand{\Dt}{\Delta t}

\newcommand{\core}{\mathrm{core}}
\newcommand{\buf}{\mathrm{buf}}

\begin{document}

\begin{frontmatter}

\title{Patch LegONet: Resolving Localized High-Frequency PDE Dynamics with Low-Frequency Operator Blocks}

\author[label1]{Jiahao Zhang\fnref{fn1}}
\ead{zhan2296@purdue.edu}

\author[label1]{Junyi Wu\fnref{fn1}}

\fntext[fn1]{Jiahao Zhang and Junyi Wu contributed equally to this work.}

\author[label1]{Yueqi Wang}

\author[label1,label2]{Guang Lin\corref{cor1}}
\cortext[cor1]{Corresponding author}
\ead{guanglin@purdue.edu}

\affiliation[label1]{organization={Purdue University, Department of Mathematics},
            city={West Lafayette},
            postcode={47906},
            state={IN},
            country={United States}}

\affiliation[label2]{organization={Purdue University, School of Mechanical Engineering},
            city={West Lafayette},
            postcode={47906},
            state={IN},
            country={United States}}

\begin{abstract}

Localized interfaces, layers, and shocks can require high spatial resolution even when most of a partial differential equation (PDE) solution remains smooth. In global spectral methods, resolving these features requires increasing the approximation order across the entire domain. LegONet builds PDE solvers from reusable mechanism blocks, but its global coefficient representation faces the same limitation. We introduce Patch LegONet, which reuses pretrained mechanism blocks on smaller physical patches to resolve localized fine scales. The solution consists of a global background and local corrections. Mapping a fixed reference representation onto smaller patches increases local resolution without enlarging the blocks. The corrections join smoothly to the global background at patch boundaries. Differential operators are rescaled analytically to account for each patch’s physical size. Nonlinear mechanisms are evaluated on the combined solution to include interactions between the background and local corrections. Patches can move, merge, or be removed during rollout while the pretrained patch blocks remain frozen. Six experiments on Allen–Cahn and viscous Burgers equations assess static, multiple, and moving patches in one to three dimensions. At comparable numbers of degrees of freedom, Patch LegONet consistently lowers the final relative solution error compared with global refinement. The reduction is approximately one order of magnitude in most tests. These results extend mechanism reuse across spatial scales and provide a practical route to local refinement in learned PDE solvers.

\end{abstract}


\end{frontmatter}

\section{Introduction}
\label{sec:introduction}

Localized fine structures can determine the resolution requirement of an entire PDE simulation. A thin interface, narrow diffusion layer, or viscous shock may occupy only a small part of the domain, while the remaining solution varies smoothly. In a global spectral representation, resolving such a feature requires increasing the approximation order across the entire domain. The resulting state is therefore sized for the most demanding local structure. This creates a mismatch between where resolution is needed and where the additional degrees of freedom are placed.

Classical numerical methods address this mismatch through local approximation. Spectral element methods place high-order spaces on individual elements \cite{Patera1984}, while adaptive mesh refinement introduces finer grids in selected regions \cite{BergerOliger1984,BergerColella1989}. Partition-of-unity and generalized finite element methods enrich a background representation with local functions \cite{MelenkBabuska1996,BabuskaMelenk1997}, and multiscale methods couple local models to a coarse solver \cite{EEngquist2003HMM}. These approaches establish a basic principle of scientific computing: spatial resolution should follow local solution complexity. For learned PDE solvers, applying this principle raises an additional question about the reuse of pretrained components.

Learned PDE solvers introduce new ways to represent solutions and approximate their evolution. Physics-informed neural networks fit solution representations using governing-equation residuals and prescribed constraints \cite{RaissiPINN2019,KarniadakisPIML2021,MoseleyFBPINNs2023,DoleanMultilevelFBPINNs2024}. DeepONet, Fourier neural operators, and related approaches learn mappings between functions for repeated prediction \cite{LuDeepONet2021,LiFNO2021,KovachkiNO2023,KovachkiSurvey2024}. Nevertheless, a solver using a fixed finite spectral state inherits the resolution limit of its retained modes. Evaluating that state on a finer grid does not add missing modes. The physical mechanism may remain unchanged while the resolution needed to represent its solution increases. A reusable learned solver should accommodate this difference without automatically requiring a new model.

LegONet provides a foundation for reuse at the level of physical mechanisms \cite{ZhangLegONet2026}. It decomposes a PDE into components such as diffusion, transport, and reaction. Each component is represented by an independently trained or analytically specified operator block. A boundary-adapted coefficient interface, called a baseplate, supplies the common state on which the blocks act. Pretrained blocks are then selected, scaled, and composed to assemble different PDE solvers. A geometry-aware extension realizes frozen mechanisms on new domains through boundary-adapted coordinates \cite{ZhangGeometryLegONet2026}. However, mechanism reuse alone does not remove the resolution limit of a fixed global baseplate. A sufficiently narrow feature can still require a larger coefficient state and blocks matched to the higher reference capacity.

Local and multiresolution operator learners offer relevant approaches to spatial reuse. Wavelet and multiwavelet architectures represent operators across several scales \cite{GuptaMWT2021,TripuraWNO2023,HeMgNO2024}. Local neural operators reuse learned time evolution on different domains \cite{LiLocalNO2024, YinHybridFEMDeepONet2022, WangNeuralFECoupling2025, YeLNOLocality2024}. Schwarz neural inference combines learned subdomain solutions through iterative interface updates and uses coordinate transformations to normalize local problems \cite{HuangDDNO2026}. These developments establish locality and coordinate transformation as useful tools for learned solvers. Instead of applying these tools to the complete PDE solution, our aim is to apply them to individual PDE mechanisms: the spatial representation should gain local resolution while the pretrained reference mechanisms remain reusable. 

We introduce Patch LegONet to make local physical resolution part of mechanism deployment. The solution is represented by a global background and a collection of local corrections,
\begin{equation}
    u_h(\mathbf{x},t)
    =
    u_G(\mathbf{x},t)
    +
    \sum_{p\in\mathcal{A}(t)}u_p(\mathbf{x},t).
    \label{eq:intro_additive_representation}
\end{equation}
Here $u_G$ represents the background, $\mathcal{A}(t)$ is the active patch set, and $u_p$ is a correction supported on patch $\Omega_p$. The global component captures the large-scale solution, while the corrections supply additional resolution near localized features. Each learned mechanism block is pretrained once on the reference patch. All physical patches reuse these same blocks with identical, fixed learned parameters, regardless of their positions or sizes. Each patch carries its own evolving coefficient state.

The resolution gain comes from changing the physical scale of the reference representation. For an isotropic patch, reference coordinates $\boldsymbol{\xi}\in[-1,1]^d$ are mapped into physical coordinates by
\begin{equation}
    \mathbf{x}
    =
    T_p(\boldsymbol{\xi})
    =
    \mathbf{c}_p+h_p\boldsymbol{\xi},
    \qquad
    K_{\mathrm{phys},p}
    \sim
    \frac{K}{h_p}.
    \label{eq:intro_scale_resolution}
\end{equation}
Here $\mathbf{c}_p$ is the patch center, $h_p$ is its half-width, and $K$ denotes the reference spectral capacity. The quantity $K_{\mathrm{phys},p}$ describes the corresponding physical frequency capacity. A smaller patch places the same reference modes over a shorter physical distance. It therefore supplies finer local resolution without increasing the reference coefficient dimension or changing the learned patch parameters. Patch size must also leave sufficient room to contain the feature and its surrounding buffer.

Three numerical interfaces connect this local representation to the PDE dynamics. First, boundary-compatible correction spaces allow the local corrections to join smoothly to the global background. A buffer region accommodates this transition. Second, analytic coordinate transformations rescale the differential mechanisms to account for each patch's physical size. Third, nonlinear mechanisms are evaluated on the combined background and correction to retain their interactions. These interfaces act outside the pretrained mechanism library, so no separate interaction model is required.

Patch management allows this local resolution to follow evolving structures. Patches may remain fixed, move with a feature, merge with nearby patches, or be removed as corrections decay. Changes in the patch configuration are handled by reconstructing the current field and projecting it onto the updated representation. The reference patch blocks remain frozen throughout these operations. We focus on local differential mechanisms and features that can be contained within buffered patches. Prescribed patches and feature-based relocation are used to assess scale reuse under controlled conditions.

The main contributions are summarized as follows.
\begin{enumerate}
    \item \textbf{Reuse across physical scales.}
    We develop a global--patch formulation that increases local physical resolution through geometric deployment. The reference patch representation and learned parameters remain fixed across patch positions and sizes.

    \item \textbf{Compatible interfaces for local mechanism evaluation.}
    We combine boundary-compatible correction spaces, analytic mechanism scaling, and nonlinear evaluation of the combined background and correction. These interfaces support local refinement within the LegONet composition framework.

    \item \textbf{Patch updates during solution evolution.}
    We provide procedures for relocating, merging, and removing patches. These operations change the local representation without retraining the reference patch library.

    \item \textbf{Controlled validation of localized resolution.}
    Six experiments on Allen--Cahn and viscous Burgers equations assess static, multiple, and moving patches in one to three dimensions. At comparable numbers of degrees of freedom, the patched solver consistently improves on global refinement, with final-error reductions of approximately one order of magnitude in most tests.
\end{enumerate}

The broader implication is that a pretrained mechanism library need not remain tied to one physical resolution. Its spatial deployment can change as localized structures evolve, while its learned patch parameters remain fixed. This extends the usefulness of mechanism-based operator libraries to simulations in which fine-scale activity is concentrated in small, changing regions. It also connects modular operator learning with the established numerical practice of allocating resolution according to the solution.

The remainder of the paper is organized as follows. Section~\ref{sec:background} reviews the LegONet coefficient interface and the setting for local scale reuse. Section~\ref{sec:methodology} develops the patch representation, mechanism transfer, correction dynamics, and deployment procedure. Section~\ref{sec:experiments} presents numerical experiments. Section~\ref{sec:discussion} discusses the results and extensions, and Section~\ref{sec:conclusion} concludes.


\section{LegONet and local scale reuse}
\label{sec:background}

LegONet builds a PDE solver from reusable blocks for individual physical mechanisms. To explain the setting for Patch LegONet, we first describe what these blocks do and how they form a solver. We then distinguish the quantities that change during a simulation from those fixed by pretraining. This distinction leads to the spatial resolution problem addressed in Section~\ref{sec:methodology}.

\subsection{A solver assembled from physical mechanisms}

Many PDEs combine several physical processes that act on the same field. For example, a reaction--diffusion equation has the form
\begin{equation}
    \partial_t u=\nu\Delta u+f(u),
    \label{eq:reaction_diffusion_example}
\end{equation}
where diffusion smooths spatial variations and the reaction $f(u)$ changes the field locally. These two processes have different mathematical actions, even though they jointly determine the evolution of one solution.

LegONet represents these actions with separate mechanism blocks \cite{ZhangLegONet2026}. A diffusion block computes the effect of diffusion on the current state. A reaction block computes the effect of the chosen reaction law. Their contributions are then assembled to advance the solution. More general equations are treated by selecting blocks for the mechanisms they contain, such as diffusion, transport, and reaction, and assigning the corresponding physical coefficients.

The reusable object is therefore the action of a mechanism on a state. A block does not predict an entire solution trajectory. It supplies the contribution needed by a numerical time integrator at the current state. This follows the familiar method-of-lines viewpoint: spatial mechanisms define a system of ordinary differential equations, and a time integrator advances that system. LegONet introduces independently constructed blocks for the individual mechanisms. A block may be learned or supplied by an analytical discrete operator.

This separation makes reuse possible at the mechanism level. For example, the same diffusion block can participate in different combinations of diffusion and reaction on a common spatial representation. Changing the diffusion strength changes its physical multiplier. The block itself continues to represent the same reference diffusion mechanism.

\subsection{A common coefficient state connects the blocks}

LegONet represents the solution by coefficients in a boundary-compatible trial space. This representation provides a common interface for the mechanism blocks. On a reference domain $\Oh$, the field is written as
\begin{equation}
    \widehat u(\xiq,t)
    =\widehat u_{\mathrm{lift}}(\xiq)+\Phihat(\xiq)a(t),
    \qquad a(t)\in\R^{r_K}.
    \label{eq:legonet_baseplate}
\end{equation}
At reference position $\xiq\in\Oh$ and time $t$, the field $\widehat u$ combines a prescribed lifting $\widehat u_{\mathrm{lift}}$ with a finite basis expansion. Choosing the approximation order $K$ determines the $r_K$ retained basis functions and the dimension of their coefficient vector $a(t)$. Collecting their values in the row vector $\Phihat(\xiq)=(\widehat\phi_1(\xiq),\ldots,\widehat\phi_{r_K}(\xiq))$ makes the product $\Phihat(\xiq)a(t)$ their weighted sum. Evaluating these rows at the numerical nodes produces the basis evaluation matrix.

The lifting supplies the prescribed nonhomogeneous boundary data, while the basis functions satisfy the corresponding homogeneous conditions. Their coefficients can therefore evolve without changing the imposed boundary values. We use a stationary lifting here, consistent with the numerical examples. To recover the coefficients from a field, a fixed linear projection $\Phat$ is applied to the remainder $\widehat u-\widehat u_{\mathrm{lift}}$. Together, the lifting, basis reconstruction, and coefficient projection form the reference baseplate. This shared interface allows all associated blocks to act on states of the same dimension and express their outputs in the same basis.

For the reaction--diffusion example, the assembled coefficient dynamics take the form
\begin{equation}
    \dot a(t)
    =\nu F_{\mathrm{diff},K}^{\theta}(a(t))
     +F_{\mathrm{react},K}^{\theta}(a(t)).
    \label{eq:legonet_coefficient_assembly}
\end{equation}
Computing the coefficient time derivative $\dot a=\mathrm{d}a/\mathrm{d}t$ requires two mechanism evaluations. The diffusion operator $F_{\mathrm{diff},K}^{\theta}$ and reaction operator $F_{\mathrm{react},K}^{\theta}$ both map the current state in $\R^{r_K}$ to a contribution in the same coefficient space. Multiplying the diffusion output by the physical diffusivity $\nu$ and adding the reaction output gives the assembled derivative. Both operators use the reference approximation order $K$, with their learned parameters collected in $\theta$.

More generally, LegONet assembles the coefficient dynamics as
\begin{equation}
    \dot a(t)
    =\sum_{i=1}^{N_{\mathrm{blk}}}
      c_iF_{i,K}^{\theta}(a(t);\mu_i).
    \label{eq:legonet_general_assembly}
\end{equation}
The sum combines $N_{\mathrm{blk}}$ selected mechanisms acting on the same coefficient state $a(t)\in\R^{r_K}$. For mechanism $i$, the block $F_{i,K}^{\theta}$ evaluates its contribution using any additional mechanism parameters collected in $\mu_i$. The scalar physical multiplier $c_i$ scales this output before assembly. As in the preceding example, $K$ identifies the reference approximation order and $\theta$ denotes the learned parameters (empty for analytical blocks). We use this notation for mechanism blocks and their assembly throughout this work. A numerical time integrator then advances $a(t)$ using the assembled contributions. Splitting methods instead advance the mechanisms in separate substeps and compose their updates \cite{Strang1968,Hairer2006}. In either case, the physical multipliers are applied outside the pretrained blocks. For later use, we also define the field action of a block,
\begin{equation}
    \mathscr{L}_{i,K}^{\theta}[v]
    :=\Phihat F_{i,K}^{\theta}(\Phat v),
    \label{eq:field_level_block}
\end{equation}
which projects a homogeneous field $v$ onto the coefficient space, applies the block and reconstructs the result.

Pretraining fits each block to the instantaneous action of its mechanism. The reference mechanism is evaluated on fields reconstructed from sampled coefficient states. Projecting each response back into coefficient space provides the target for fitting the block parameters. Each sample therefore supplies an operator evaluation, without requiring a full solution trajectory \cite{ZhangLegONet2026}.

Structure preservation is a central feature of the LegONet block design. To retain prescribed dissipation or conservation, learned scalar generators are combined with fixed structure operators. Representative dissipative E-blocks and Hamiltonian H-blocks \cite{GreydanusHNN2019} have the forms
\begin{equation}
    F_{E,K}^{\theta}(a)
    =-G\nabla_a E^{\theta}(a),
    \qquad
    F_{H,K}^{\theta}(a)
    =J\nabla_a H^{\theta}(a).
    \label{eq:legonet_structured_blocks}
\end{equation}
The learned generators $E^{\theta},H^{\theta}:\R^{r_K}\to\R$ assign scalar energy and Hamiltonian values to a coefficient state. Differentiating them with respect to $a$, as indicated by $\nabla_a$, gives vectors in coefficient space. Fixed matrices $G,J\in\R^{r_K\times r_K}$ convert these gradients into the dissipative operator $F_{E,K}^{\theta}$ and Hamiltonian operator $F_{H,K}^{\theta}$. Symmetric positive-semidefinite structure is imposed through $G=G^{\mathsf T}\succeq0$, while $J^{\mathsf T}=-J$ imposes skew symmetry. These matrices include any mass-matrix factors required by the coefficient representation. Along the isolated continuous-time E-block flow, $E^{\theta}$ is nonincreasing. Along the corresponding H-block flow, $H^{\theta}$ is conserved. These properties follow from the block construction and hold throughout the entire rollout process. Detailed block architectures are given in the original LegONet formulation \cite{ZhangLegONet2026}.

\subsection{Mechanism reuse and the resolution limit}

After pretraining, LegONet advances a solution by evolving the coefficient state while keeping the learned mechanism parameters fixed. The same blocks can act on different admissible states and be combined into solvers that share these mechanisms. This reuse takes place within the coefficient interface defined by the reference baseplate.

A fixed reference baseplate also imposes a spatial resolution limit. With the approximation order $K$ fixed, updating the $r_K$ coefficients changes the amplitudes of the retained basis functions but does not introduce additional modes. Evaluating the reconstructed field at more points only samples that representation more densely. Consequently, even an accurately learned mechanism cannot recover spatial structures that the chosen approximation space cannot represent.

This limitation becomes important when fine structures occupy only a small part of the domain. A smooth background may require relatively few coefficients, while a thin interface or localized layer requires much finer resolution. Increasing the global approximation order allocates additional degrees of freedom throughout the domain, although the extra resolution is needed only locally. The larger coefficient state also requires mechanism blocks compatible with the enlarged coefficient space.

The resulting challenge is to increase local physical resolution while retaining the pretrained reference library. Patch LegONet addresses this challenge by deploying the same reference blocks on smaller physical patches, as developed in Section~\ref{sec:methodology}.


\section{Patch LegONet methodology}
\label{sec:methodology}

\begin{figure}[!t]
\centering
\includegraphics[width=0.96\linewidth]{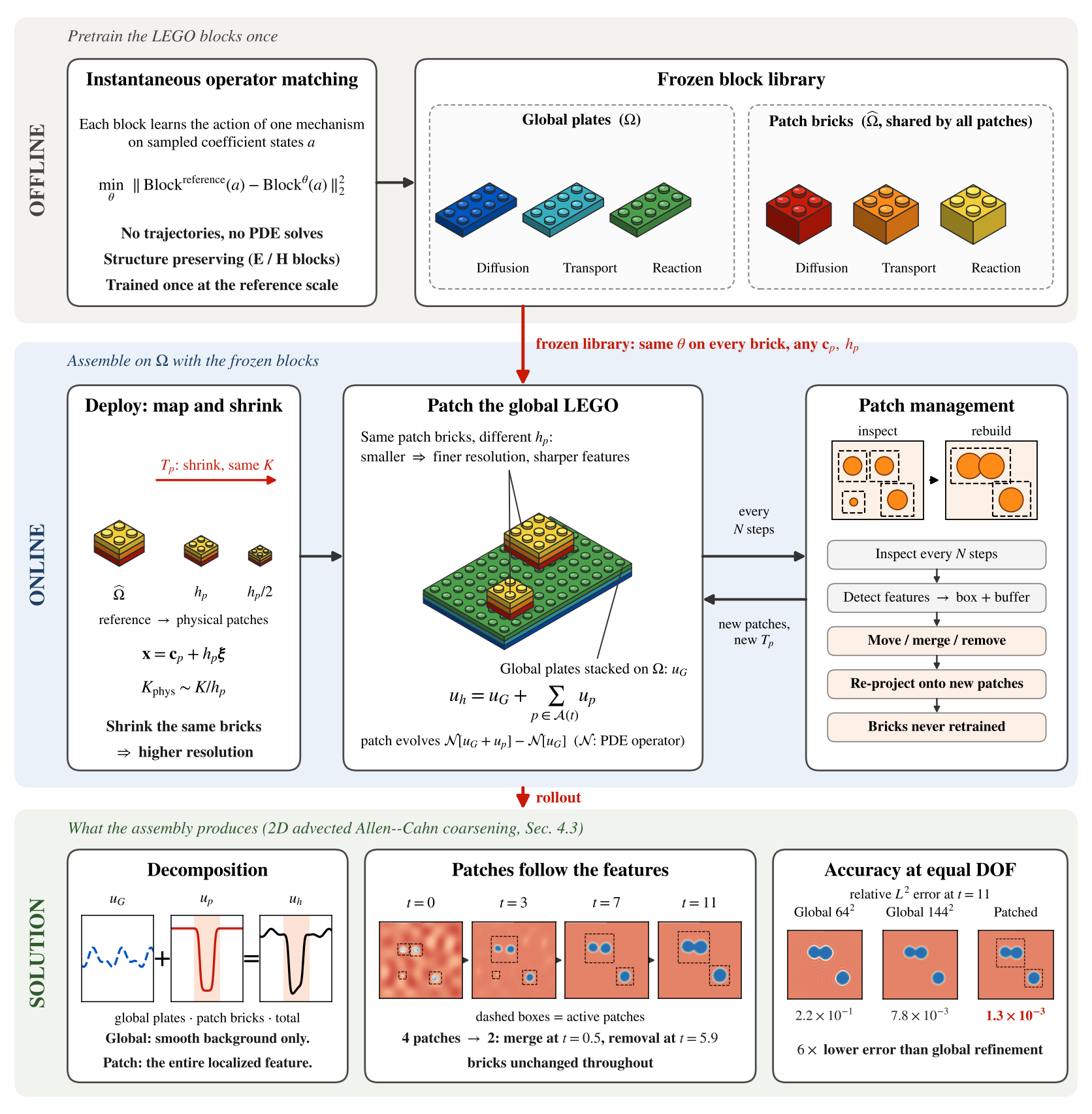}
\caption{\textbf{Patch LegONet workflow.}
\emph{Offline}: mechanism blocks are fitted once on the reference baseplates by instantaneous operator matching, giving a frozen library of global plates (on $\Omega$) and patch bricks (on $\widehat{\Omega}$, shared by all patches).
\emph{Online}: the patch bricks are mapped to physical patches by $T_p$ and shrunk to the required local resolution $K_{\mathrm{phys}}\sim K/h_p$; the global plates are stacked on $\Omega$ and the patch bricks are stacked locally on top, so that $u_h=u_G+\sum_p u_p$ with each patch evolving the operator difference $\mathcal N[u_G+u_p]-\mathcal N[u_G]$; patches are inspected, moved, merged, or removed by projection while the blocks remain frozen.
\emph{Solution}: the global state carries only the smooth background and each patch carries the entire localized feature; for the 2D advected Allen--Cahn coarsening of Section~\ref{subsec:exp_moving}, the patches follow, merge, and vanish with the droplets, and the final error is six times lower than that of global refinement at a comparable coefficient count.}
\label{fig:workflow}
\end{figure}

Patch LegONet separates the reference representation of physical mechanisms from the physical scale at which they are used. The solution combines a global background with local corrections carried by physical patches. Each patch has its own coefficient state, while the reference patch basis and pretrained mechanism parameters remain fixed. We first define the global--patch representation and its resolution gain. We then develop boundary compatibility, the scaling of the reference operators, the nonlinear coupling, and the deployment and update of patches.
Figure~\ref{fig:workflow} summarizes the complete workflow.

\subsection{Global background and local corrections}
\label{subsec:global_patch_representation}

The construction starts from a spatial decomposition that is common to all mechanisms in the PDE. On a physical domain $\Omega\subset\R^d$, with $d$ spatial dimensions, we approximate the solution by
\begin{equation}
    u_h(\x,t)
    =
    u_G(\x,t)
    +
    \sum_{p\in\Acal(t)}u_p(\x,t).
    \label{eq:additive_global_patch}
\end{equation}
At physical position $\x$ and time $t$, the approximate field $u_h$ combines the global background $u_G$ with corrections indexed by the active patch set $\Acal(t)$. Each correction $u_p$ is supported on a physical patch $\Omega_p\subset\Omega$. The background remains present throughout the domain, including inside the patches, while the corrections supply additional local approximation capacity.

We represent the background on a global baseplate as
\begin{equation}
    u_G(\x,t)
    =
    u_{\mathrm{lift},G}(\x)
    +
    \Phi_G(\x)a_G(t),
    \qquad
    a_G(t)\in\R^{r_G}.
    \label{eq:global_component}
\end{equation}
The row $\Phi_G(\x)$ collects the values of $r_G$ global basis functions, whose evolving coefficients form the vector $a_G(t)$. The stationary lifting $u_{\mathrm{lift},G}$ supplies the prescribed nonhomogeneous boundary data, and the basis functions satisfy the corresponding homogeneous conditions.

All local patches are constructed from one reference domain $\Oh=[-1,1]^d$. For the isotropic patches considered here, the reference-to-physical map takes the form
\begin{equation}
    \Omega_p=T_p(\Oh),
    \qquad
    T_p(\xiq)=\cp+h_p\xiq.
    \label{eq:patch_map_general}
\end{equation}
The map $T_p$ sends a reference coordinate $\xiq\in\Oh$ to the physical patch $\Omega_p$. Its center $\cp\in\R^d$ determines the location, and its positive half-width $h_p$ determines the size of the patch. Patch positions and sizes are held fixed between the update events described in Section~\ref{subsec:patch_management}.

Using the same reference basis on every patch, we define
\begin{equation}
    u_p(\x,t)=
    \begin{cases}
        \Phihat\!\left(T_p^{-1}(\x)\right)a_p(t),
        & \x\in\Omega_p,\\
        0,
        & \x\in\Omega\setminus\Omega_p.
    \end{cases}
    \label{eq:patch_component}
\end{equation}
The inverse coordinate $T_p^{-1}(\x)=(\x-\cp)/h_p$ identifies where the reference basis row $\Phihat$ is evaluated. At the fixed approximation order $K$, its $r_K$ basis functions reconstruct the correction from the coefficients in $a_p(t)\in\R^{r_K}$. Patch corrections have zero lifting. The coefficient constraints required for a smooth zero extension are given in Section~\ref{subsec:patch_compatibility}.

The global and patch bases may differ, and each baseplate uses its own compatible mechanism blocks. For the reference patch baseplate we fix one library
$\Bcal_K^\theta=\{F_{i,K}^\theta\}_{i=1}^{N_{\mathrm{blk}}}$
in the notation of Section~\ref{sec:background}. Each block is pretrained once on this baseplate, and all physical patches share its learned parameters and basis while carrying separate coefficient states.

\subsection{Local physical resolution at fixed reference capacity}
\label{subsec:local_resolution}

The patch map changes the physical scale of the reference approximation space. Reducing the half-width $h_p$ places these basis functions on a smaller physical region. Their variations therefore occur over shorter distances, allowing the same reference representation to describe finer local structure.

A single reference mode illustrates this geometric effect. In one spatial dimension, the map $x=c_p+h_p\xi$ transforms an illustrative cosine mode according to
\begin{equation}
    \widehat\phi_k(\xi)=\cos(\pi k\xi),
    \qquad
    \phi_{p,k}(x)
    =
    \cos\!\left(
        \pi k\frac{x-c_p}{h_p}
    \right).
    \label{eq:reference_mode_on_patch}
\end{equation}
The integer $k=1,\ldots,K$ indexes the retained reference modes, and $\phi_{p,k}$ denotes the mode expressed in physical coordinates. Translating the patch center $c_p$ changes its location, while the half-width $h_p$ sets its physical wavelength $2h_p/k$. Equivalently, its physical angular wavenumber is $\pi k/h_p$. Thus the highest retained reference modes provide the effective physical frequency capacity
\begin{equation}
    K_{\mathrm{phys},p}
    \sim
    \frac{K}{h_p}.
    \label{eq:scale_frequency_method}
\end{equation}
The same scaling holds for polynomial baseplates. Mapping a reference polynomial onto a patch preserves its degree and multiplies every node spacing by $h_p$.

This relation also guides the choice of patch size. To resolve a target characteristic physical frequency $M_p$, we require
\begin{equation}
    K_{\mathrm{phys},p}\gtrsim M_p,
    \qquad\text{which gives}\qquad
    h_p\lesssim\frac{K}{M_p}.
    \label{eq:h_from_frequency}
\end{equation}

A lower bound on $h_p$ follows from the buffer requirement in Section~\ref{subsec:patch_compatibility}.

\subsection{Boundary-compatible correction spaces}
\label{subsec:patch_compatibility}

The local correction in \eqref{eq:patch_component} is extended
by zero outside its patch. Its boundary values and derivatives
must therefore match the zero exterior field to the required
order. These homogeneous boundary constraints apply only to
the correction. The global background can remain nonzero at
the artificial patch boundary.

The required compatibility depends on the differential
mechanisms in the PDE. Let $q\ge 1$ denote the highest spatial
derivative order, and let $\ell$ denote the continuity order
imposed on the extended correction. For a correction that is
smooth up to the boundary from inside the patch, the condition
\begin{equation}
    \ell\ge q-1
    \label{eq:patch_compatibility_order}
\end{equation}
is sufficient to prevent singular contributions supported
on the patch boundary when taking derivatives through order
$q$. This condition allows an ordinary jump in the
highest-order derivatives. If those derivatives must also
be continuous across the boundary, we require $\ell\ge q$
instead.

We impose the boundary constraints through the reference
approximation space. The compatible coefficient space is
\begin{equation}
    \mathcal Z_K^{(\ell)}
    :=
    \left\{
        a\in\R^{r_K}:
        D_{\xi}^{\alpha}
        \big[\Phihat(\xiq)a\big]=0
        \quad\text{on }\partial\Oh,
        \quad |\alpha|\le\ell
    \right\}.
    \label{eq:compatible_patch_space}
\end{equation}
Here $\alpha$ is a multi-index of total order $|\alpha|$.
The conditions are linear in $a$, so
$\mathcal Z_K^{(\ell)}$ is a linear subspace of
$\R^{r_K}$.

\begin{proposition}[Compatible zero extension]
\label{prop:smooth_zero_extension}
Suppose $\Phihat(\cdot)a_p\in C^{\ell+1}(\overline{\Oh})$
and $a_p\in\mathcal Z_K^{(\ell)}$. Then the correction
defined by \eqref{eq:patch_component} belongs to
$C^\ell(\Omega)$. If $\ell\ge q-1$, its weak spatial
derivatives through order $q$ coincide almost everywhere
with the corresponding interior derivatives extended
by zero.
\end{proposition}

\begin{proof}
The affine patch map preserves the vanishing boundary
conditions. Hence the interior values and derivatives
through order $\ell$ match the zero exterior values,
giving a $C^\ell$ extension. To compute a weak derivative
of order at most $q$, repeated integration by parts
produces boundary terms involving derivatives of the
correction of order at most $q-1$. These terms vanish
because $q-1\le\ell$, leaving only the corresponding
interior derivative extended by zero.
\end{proof}

The numerical examples in this work contain second-order
diffusion, so we choose $\ell=1$. We satisfy the required
boundary constraints directly by using the clamped
Shen--Legendre basis \cite{Shen1994, ShenTangWang2011}. In one dimension,
its basis functions can be written, up to a fixed
normalization, as
\begin{equation}
    \psi_k(\xi)
    =
    L_k(\xi)
    -
    \frac{2(2k+5)}{2k+7}L_{k+2}(\xi)
    +
    \frac{2k+3}{2k+7}L_{k+4}(\xi),
    \qquad k=0,\ldots,K-1.
    \label{eq:clamped_shen_basis}
\end{equation}
Here $L_j$ denotes the Legendre polynomial of degree $j$,
and the index range retains $K$ one-dimensional basis
functions. This combination gives
$\psi_k(\pm1)=\psi_k'(\pm1)=0$, where the prime denotes
differentiation with respect to $\xi$. After mapping
onto a physical patch, any linear combination therefore
extends by zero as a $C^1$ function.

For a $d$-dimensional reference patch, we use tensor
products of these functions:
\begin{equation}
    \widehat\phi_{\bm{k}}(\xiq)
    =
    \prod_{j=1}^{d}\psi_{k_j}(\xi_j),
    \qquad
    \bm{k}\in\{0,\ldots,K-1\}^{d}.
    \label{eq:clamped_shen_tensor_basis}
\end{equation}
The index tuple $\bm{k}=(k_1,\ldots,k_d)$ selects one
basis function in each coordinate direction. On every
face of the reference patch, the corresponding factor
and its first derivative vanish. Consequently, each
tensor-product function and all its first partial
derivatives vanish on the boundary.

Collecting these $r_K=K^d$ functions in the reference
basis row $\Phihat$ gives the reconstruction used in
\eqref{eq:patch_component}. Every coefficient vector
$a_p\in\R^{r_K}$ then satisfies the required compatibility,
so $\mathcal Z_K^{(1)}=\R^{r_K}$. The patch state can
therefore be evolved directly in this basis. The fixed
projection $\Phat$ introduced in
Section~\ref{sec:background} also uses this basis, so
projecting a field contribution returns coefficients
whose reconstruction satisfies the same boundary
constraints.

The patch must also leave room for the correction to
decay to zero. We divide the reference patch into a core
and a surrounding buffer,
\begin{equation}
    \Oh_{\core}=[-\rho,\rho]^d,
    \qquad
    \Oh_{\buf}=\Oh\setminus\Oh_{\core},
    \qquad 0<\rho<1.
    \label{eq:reference_core_buffer}
\end{equation}
The parameter $\rho$ sets the core half-width in each
reference coordinate. The physical core
$\Omega_{p,\core}=T_p(\Oh_{\core})$ contains the
localized feature. The surrounding buffer
$\Omega_{p,\buf}=T_p(\Oh_{\buf})$ provides room for
the correction and its required derivatives to
return to zero at the patch boundary.

Let $r_p$ denote a containing half-width of the
feature about $\cp$, measured in the maximum norm
over the interval until the next patch update.
Containment within the core requires
\begin{equation}
    r_p\le\rho h_p.
    \label{eq:buffer_size_constraint}
\end{equation}
Combining this condition with the resolution
criterion \eqref{eq:h_from_frequency} gives
\begin{equation}
    \frac{r_p}{\rho}
    \le h_p
    \lesssim
    \frac{K}{M_p}.
    \label{eq:admissible_patch_interval}
\end{equation}
The lower bound leaves room for the feature and its buffer. The upper bound provides the required resolution at fixed $K$.
For a feature of half-extent $r_p$, we take the smallest admissible patch,
$h_p=r_p+b$, which maximizes $K/h_p$; hence $\rho=1-b/h_p$. The padding $b$
is chosen so that the decaying tail of the feature falls below a small
tolerance before the patch boundary.

Other reference bases can be used if they satisfy
\eqref{eq:compatible_patch_space} at the required
order $\ell$, with the mechanism blocks and projection
built on the same baseplate.

\subsection{Scaling of the reference operators}
\label{subsec:analytic_pullback}

The reference operators are defined on $\Oh$. On a physical patch $\Omega_p$, the spatial derivatives are rescaled by the half-width $h_p$. The resulting factors multiply the reference operators, whose pretrained parameters remain fixed.

Consider a physical field $w(\x,t)$ on $\Omega_p$ and
its reference representation
$\widehat w(\xiq,t)=w(T_p(\xiq),t)$. While the patch
center and half-width remain fixed, the chain rule gives
\begin{equation}
    \big(D_x^\alpha w\big)(T_p(\xiq),t)
    =
    h_p^{-|\alpha|}
    D_\xi^\alpha\widehat w(\xiq,t).
    \label{eq:general_patch_derivative}
\end{equation}
The operators
$D_x^\alpha$ and $D_\xi^\alpha$ differentiate in the
physical and reference coordinates, respectively.
Each spatial derivative therefore introduces one factor
of $h_p^{-1}$. Time differentiation is unchanged because
the patch geometry is fixed between updates.

For a mechanism whose coordinate transformation
produces a single scalar factor, we define the
corresponding block on patch $p$ by
\begin{equation}
    F_{i,p}^{\theta}(a;\mu_i)
    :=
    h_p^{-m_i}
    F_{i,K}^{\theta}(a;\mu_i).
    \label{eq:primitive_block_scaling}
\end{equation}
The block $F_{i,p}^{\theta}$ applies the reference
operator $F_{i,K}^{\theta}$ to the coefficient state
$a$, using the same learned parameters $\theta$ and
mechanism parameters $\mu_i$. The exponent $m_i$
records the power of $h_p^{-1}$ generated by the
coordinate transformation. For diffusion $\Delta u$,
$m_i=2$; for a first-derivative transport term,
$m_i=1$; and for a pointwise reaction, $m_i=0$.
Thus the physical multiplier $c_i$ in
\eqref{eq:legonet_general_assembly} becomes
$c_i h_p^{-m_i}$ on patch $p$.

Products of spatial derivatives acquire the product of their factors. Terms with different powers of $h_p$ are rescaled separately before assembly, following the mechanism decomposition of Section~\ref{sec:background}.
In patch coordinates, the reaction--diffusion equation
\eqref{eq:reaction_diffusion_example} becomes
\begin{equation}
    \partial_t\widehat u
    =
    \frac{\nu}{h_p^2}\Delta_\xi\widehat u
    +
    f(\widehat u).
    \label{eq:reaction_diffusion_pullback}
\end{equation}
Here $\Delta_\xi$ is the reference Laplacian. The physical
diffusivity $\nu$ is multiplied by $h_p^{-2}$,
while the reaction law $f$ retains its original form.
The same reference diffusion and reaction blocks can
therefore be used with the modified multipliers.
Similarly, the Burgers transport term
$u\,\partial_x u$ receives the factor $h_p^{-1}$.

Smaller patches make the equations stiffer. For a fixed
reference discretization, explicit diffusion stability
depends on $\nu\Dt/h_p^2$, where $\Dt$ is the physical time step. The time step must therefore account for the patch size.

\subsection{Coupling the global background and local corrections}
\label{subsec:total_jet_assembly}

The coordinate factors in
Section~\ref{subsec:analytic_pullback} determine how
the reference mechanisms are rescaled on a physical
patch. Their application to a local correction must
also account for the global background. For a nonlinear
mechanism, evaluating the background and correction
separately generally omits their interaction. We formulate the coupling for patches with nonoverlapping interiors. On patch $p$, the total field then consists of the background and its local
correction. In reference coordinates, these fields are
\begin{equation}
    \widehat u_{G,p}(\xiq,t)
    =
    u_G(T_p(\xiq),t),
    \qquad
    \widehat v_p(\xiq,t)
    =
    \Phihat(\xiq)a_p(t).
    \label{eq:local_background_and_correction}
\end{equation}
The background $\widehat u_{G,p}$ is evaluated from
the global representation and includes the lifting
in \eqref{eq:global_component}. The reference basis
$\Phihat$ reconstructs the correction from its patch
coefficients. Because the two bases may differ, we
combine their field values instead of coefficients.

Let $\widehat{\Ncal}_p$ denote the state-dependent
spatial operator in patch coordinates, including the
scale factors from
Section~\ref{subsec:analytic_pullback}. We define
the operator difference assigned to the patch by
\begin{equation}
    \widehat{\mathcal D}_p
    \big[\widehat u_{G,p},\widehat v_p\big]
    :=
    \widehat{\Ncal}_p
    \big[\widehat u_{G,p}+\widehat v_p\big]
    -
    \widehat{\Ncal}_p
    \big[\widehat u_{G,p}\big].
    \label{eq:patch_unresolved_dynamics}
\end{equation}
Both evaluations use the same spatial operator.
Subtracting the background operator value retains
the interaction with the correction. The difference
vanishes when the correction is zero.

We implement this rule through a coupled coefficient
block. It uses the current patch coefficients and
the background values on the patch, and returns a
vector in the correction coefficient space.
At $Q$ fixed reference quadrature nodes
$\{\xiq_q\}_{q=1}^{Q}$, the required samples are
\begin{equation}
    g_q=\widehat u_{G,p}(\xiq_q,t),
    \qquad
    v_q=(Ba_p)_q,
    \qquad
    B_{qk}=\widehat\phi_k(\xiq_q).
    \label{eq:patch_nodal_inputs}
\end{equation}
The basis evaluation matrix
$B\in\R^{Q\times r_K}$ reconstructs the correction
at these nodes. The vectors $g,v\in\R^Q$ contain
the background and correction samples, respectively.
The background retains its global representation
and its prescribed lifting.

For a mechanism evaluated from field samples, let
$\widehat{\mathcal Q}_{i,K}^{\theta}:
\R^Q\to\R^{r_K}$ denote its reference coefficient
operator. It combines the analytical or learned
mechanism with the fixed spatial operations and
projection associated with the reference baseplate.
With the mechanism parameters $\mu_i$ suppressed,
the reference patch block with zero lifting satisfies
$F_{i,K}^{\theta}(a)
=\widehat{\mathcal Q}_{i,K}^{\theta}(Ba)$.
The coupled block applies this same operator to
the total and background samples:
\begin{equation}
    \Ccal_{i,p}(a_G,a_p)
    =
    h_p^{-m_i}
    \left[
        \widehat{\mathcal Q}_{i,K}^{\theta}(g+Ba_p)
        -
        \widehat{\mathcal Q}_{i,K}^{\theta}(g)
    \right].
    \label{eq:coupled_nodal_block}
\end{equation}
Here $m_i$ is the coordinate-scaling exponent defined
in Section~\ref{subsec:analytic_pullback}, and the
dependence on $a_G$ enters through $g$.
Both evaluations use identical learned parameters
and fixed reference operations. The coupling therefore
introduces no additional trainable parameters.
When $a_p=0$, the two evaluations coincide and
$\Ccal_{i,p}=0$.

The reaction and transport operators used in the
numerical experiments illustrate this construction.
For the analytical Allen--Cahn reaction
$f(u)=u-u^3$ in
Sections~\ref{subsec:exp_ac1d}
and~\ref{subsec:exp_ac2d3d}, the coupled block is
\begin{equation}
    \Ccal_{\mathrm{react},p}(a_G,a_p)
    =
    \Phat\!\left[
        f(\widehat u_{G,p}+\widehat v_p)
        -
        f(\widehat u_{G,p})
    \right].
    \label{eq:patch_reaction_block}
\end{equation}
The prescribed reaction is evaluated at the patch
nodes, and the fixed reference projection $\Phat$
expresses the difference in the compatible patch
basis. This analytical mechanism requires no training.
Its dependence on the background is explicit:
\begin{align}
    f(\widehat u_{G,p}+\widehat v_p)
    -f(\widehat u_{G,p})
    &=
    (1-3\widehat u_{G,p}^{\,2})\widehat v_p
    -3\widehat u_{G,p}\widehat v_p^{\,2}
    -\widehat v_p^{\,3}.
    \label{eq:reaction_cross_terms}
\end{align}
The terms containing $\widehat u_{G,p}$ describe
the interaction between the background and correction.

The Burgers experiment in
Section~\ref{subsec:exp_burgers1d} uses a pretrained
transport mechanism. The LegONet block contains a
learned scalar function
$\mathfrak h_\theta:\R\to\R$.
Its derivative
$N_{\mathrm B}^{\theta}
=\mathrm d\mathfrak h_\theta/\mathrm du$
approximates the flux $u^2/2$
\cite{ZhangLegONet2026}.
For the right-hand-side transport term $-u\,u_x$,
the coupled coefficient operator is
\begin{equation}
    \Ccal_{\mathrm{adv},p}(a_G,a_p)
    =
    -\frac{1}{h_p}\Phat\!\left[
        \partial_\xi
        \left(
            N_{\mathrm B}^{\theta}
            (\widehat u_{G,p}+\widehat v_p)
            -
            N_{\mathrm B}^{\theta}
            (\widehat u_{G,p})
        \right)
    \right].
    \label{eq:patch_block}
\end{equation}
The factor $h_p^{-1}$ follows from the coordinate
scaling of the first derivative. The weak formulation used to compute
this coefficient vector is given in
\ref{app:patch_nodal_realization}.

For a linear mechanism, the background terms cancel.
The coupled block then reduces to
\begin{equation}
    \Ccal_{i,p}(a_G,a_p;\mu_i)
    =
    h_p^{-m_i}F_{i,K}^{\theta}(a_p;\mu_i).
    \label{eq:patch_linear_block}
\end{equation}
In the reaction--diffusion example of
Section~\ref{sec:background}, the diffusion difference
is therefore
$\nu h_p^{-2}\Delta_\xi\widehat v_p$.
The corresponding coefficient implementation is
specified in \ref{app:patch_diffusion}.

Each coupled block returns a vector in $\R^{r_K}$,
so its output can be summed directly in coefficient
space. The assembled patch equation is
\begin{equation}
    \dot a_p
    =
    \sum_{i=1}^{N_{\mathrm{blk}}}
        c_i\Ccal_{i,p}(a_G,a_p;\mu_i)
    +
    \Phat[\widehat s_p].
    \label{eq:patch_coefficient_dynamics}
\end{equation}
The physical multipliers $c_i$ and mechanism parameters
$\mu_i$ have the meanings introduced in
\eqref{eq:legonet_general_assembly}.
The prescribed source $\widehat s_p$ is assigned
to patch $p$ and expressed in reference coordinates.
It is added separately because a state-independent
term would cancel in
\eqref{eq:patch_unresolved_dynamics}.
The block outputs and projected source are expressed
in the compatible basis, so the correction retains
the boundary conditions of
Section~\ref{subsec:patch_compatibility}.

The assembled experimental operators are given in
\ref{app:patch_experimental_operators}.
Other local nonlinear mechanisms follow the same
coupling rule when their reference implementation
supports evaluation on the combined field and the
background. The following subsection describes
coefficient initialization and the transfer of states
when the patch configuration changes.

\subsection{Patch deployment and updates}
\label{subsec:patch_management}
\label{subsec:summary_methodology}

The preceding subsections specify the representation
and its coefficient operators. We now describe how
to initialize the solver and update the patch
configuration as features evolve. These operations
change the geometry and coefficient states while the pretrained reference blocks remain fixed.

Initialization starts from the known field $u_0$.
Let $u_{\mathrm{bg}}^0$ denote the selected smooth
background, including the prescribed boundary lifting. Using the fixed coefficient projection $P_G$
associated with the global basis, we set
$a_G^0=P_G[u_{\mathrm{bg}}^0-u_{\mathrm{lift},G}]$
and reconstruct
$u_G^0=u_{\mathrm{lift},G}+\Phi_Ga_G^0$. The initial patches may be prescribed or selected from $u_0$ by the detection procedure below.
Their final geometry must satisfy the requirements
of Section~\ref{subsec:patch_compatibility} and
the nonoverlap condition of
Section~\ref{subsec:total_jet_assembly}.
For each patch, we form the initial difference
\begin{equation}
    \widehat d_p^0(\xiq)
    =
    u_0(T_p^0(\xiq))
    -
    u_G^0(T_p^0(\xiq)),
    \label{eq:patch_initial_defect}
\end{equation}
where $T_p^0$ is the initial patch map.
The corresponding correction coefficients are
\begin{equation}
    a_p^0
    =
    \Phat[\widehat d_p^0].
    \label{eq:patch_initialization}
\end{equation}

Between patch updates, the active set and patch
geometry remain fixed. The global coefficients follow
\eqref{eq:legonet_general_assembly} on the global
baseplate, while the patch coefficients follow
\eqref{eq:patch_coefficient_dynamics}.
The global state advances independently and supplies
the background values for the patch operators. Let $t_n=n\Dt$, where $\Dt$ is the physical time step. In the experimental scheme, the global state is advanced first. Its values at $t_n$ and $t_{n+1}$
are then evaluated at the current patch nodes and
used to advance the patch states. The time integrator and background sampling are specified in \ref{app:patch_time_integration}.

For evolving features, we inspect the patch
configuration at prescribed times.
A configurable feature monitor is
\begin{equation}
    A(\x,t)
    =
    \max_{|\alpha|\le M}
    \frac{
        \Delta x_G(\x)^{|\alpha|}
        \left|
            D_x^\alpha
            \big(u_h(\x,t)-u_\infty\big)
        \right|
    }{U}.
    \label{eq:truncation_monitor}
\end{equation}
The largest derivative order included is $M$,
and the multi-index $\alpha$ specifies the spatial
derivatives. The local spacing scale
$\Delta x_G(\x)$ and amplitude scale $U>0$
normalize these quantities.
The prescribed value $u_\infty$ is the far-field
state. For $M=0$, the monitor measures the
deviation from that state. Thresholding the monitor identifies feature regions in the current numerical field. A one-sided threshold can be used when a particular phase is to be selected. The moving Allen--Cahn
experiment in Section~\ref{subsec:exp_moving} uses this zeroth-order phase detection. Detection uses $u_0$ at initialization and $u_h$
at later inspections.

For each connected feature region, we form a
bounding box with coordinate intervals
$[\ell_j,r_j]$, $j=1,\ldots,d$.
An initial isotropic candidate patch is defined by
\begin{equation}
    c_{p,j}
    =
    \frac{\ell_j+r_j}{2},
    \qquad
    h_p
    =
    \max\!\left\{
        h_{\min},
        \frac12\max_j(r_j-\ell_j)+b
    \right\}.
    \label{eq:patch_box_geometry}
\end{equation}
The prescribed physical padding $b>0$ leaves room
around the detected feature, while $h_{\min}>0$
sets the minimum half-width.
The final center and half-width must keep the patch
inside the physical domain and satisfy the core
containment condition
\eqref{eq:buffer_size_constraint}.
These requirements are checked after geometric
adjustments. The containment and resolution conditions in
\eqref{eq:admissible_patch_interval} guide the choice
of patch sizes and reference capacity.
The detector supplies feature locations and extents.
A target frequency $M_p$, when specified, is a
separate resolution requirement.

To satisfy the coupling assumption in
Section~\ref{subsec:total_jet_assembly}, the final
physical patches must have nonoverlapping interiors.
Intersecting candidates, or candidates within a
prescribed merging gap, are combined.
Their geometry is recomputed, and the checks are
repeated after any size adjustment.

Rebuilding is triggered when the patch count changes or when a matched patch's center or half-width changes beyond the prescribed geometric tolerance. The inspection interval, phase threshold, and removal tolerances used in the moving experiment are given in \ref{app:exp_moving}. A change in configuration requires a transfer
of the current coefficient states. Use superscripts $-$ and $+$ for the representations before and after the transfer. The old total field
$u_h^-$ and background $u_G^-$ can therefore be
evaluated at the new patch nodes.
For every patch in the new active set $\Acal^+$,
we compute
\begin{equation}
    a_p^+
    =
    \Phat\!\left[
        \big(u_h^- - u_G^-\big)\circ T_p^+
    \right],
    \qquad p\in\Acal^+.
    \label{eq:patch_rebuild_local}
\end{equation}
The new zero-extended corrections
$u_p^+$ are reconstructed by using \eqref{eq:patch_component}.
The remaining field is then encoded on the global
baseplate:
\begin{equation}
    a_G^+
    =
    P_G\!\left[
        u_h^-
        -
        \sum_{p\in\Acal^+}u_p^+
        -
        u_{\mathrm{lift},G}
    \right].
    \label{eq:patch_rebuild_global}
\end{equation}
Subtracting the stationary lifting gives the
homogeneous field represented by the global
coefficients.
This rebuild redistributes the current solution
between the background and local corrections,
up to projection error. It can therefore change both global and patch
coefficient states. Relocation, resizing, merging, and removal use
the same transfer procedure.

Algorithm~\ref{alg:patch_legonet} summarizes the
complete solver. Here $N_T$ is the number of time
steps, with final time $T=N_T\Dt$.

\begin{algorithm}[H]
\caption{Patch LegONet solver}
\label{alg:patch_legonet}
\begin{algorithmic}[1]
\Require Initial field $u_0$, selected full background
$u_{\mathrm{bg}}^0$, global and reference patch
baseplates, frozen mechanism libraries, physical
parameters, $\Dt$, $N_T$, and patch selection and
inspection rules.
\State Compute
$a_G^0=P_G[u_{\mathrm{bg}}^0-u_{\mathrm{lift},G}]$
and reconstruct $u_G^0$.
\State Select prescribed or detected initial patches;
finalize their geometry and enforce the containment
and nonoverlap requirements.
\State Initialize $a_p^0$ using
\eqref{eq:patch_initial_defect}--\eqref{eq:patch_initialization}.
\If{the initial patches are detected}
    \State Refit $a_G^0$ using
    \eqref{eq:patch_rebuild_global}, with
    $u_h^-=u_0$ and $u_p^+=u_p^0$.
\EndIf
\For{$n=0,\ldots,N_T-1$}
    \State Advance the global state to $a_G^{n+1}$
    with the current patch geometry fixed.
    \State Evaluate $u_G^n$ and $u_G^{n+1}$
    at the current physical patch nodes.
    \State Advance each patch to $a_p^{n+1}$
    using \eqref{eq:patch_coefficient_dynamics}
    and the sampled background.
    \If{an inspection is scheduled at $t_{n+1}$}
        \State Reconstruct $u_h$ and identify
        feature regions from the current field.
        \State Form candidate patches using
        \eqref{eq:patch_box_geometry}.
        \State Retain patches whose features have
        disappeared until their corrections satisfy
        the removal tolerance.
        \State Finalize the geometry of all candidates,
        including retained patches; merge and recheck
        until the containment and nonoverlap
        requirements are satisfied.
        \If{the patch count changes or geometric
        tolerances are exceeded}
            \State Retain the current geometry,
            coefficients, and fields $u_h^-$ and $u_G^-$.
            \State Construct the new active set
            $\Acal^+$ and maps $T_p^+$.
            \For{$p\in\Acal^+$}
                \State Compute $a_p^+$ using
                \eqref{eq:patch_rebuild_local}.
            \EndFor
            \State Reconstruct the new corrections
            and compute $a_G^+$ using
            \eqref{eq:patch_rebuild_global}.
            \State Replace the states and geometry
            at $t_{n+1}$ and refresh the
            geometry-dependent evaluation data.
        \EndIf
    \EndIf
\EndFor
\State \Return The final field
$u_h=u_G+\sum_{p\in\Acal(T)}u_p$.
\end{algorithmic}
\end{algorithm}


\section{Numerical experiments}
\label{sec:experiments}

In order to evaluate whether the fixed reference patch blocks
can resolve localized fine structures through their
physical deployment, six experiments on Allen--Cahn
and viscous Burgers equations examine local resolution,
reuse of a learned nonlinear mechanism, and changes
in the patch configuration during evolution. Section~\ref{subsec:exp_single} considers one static
patch in one to three dimensions.
Section~\ref{subsec:exp_multi} examines several
localized features, each represented by a separate
patch state.
Section~\ref{subsec:exp_moving} then considers patches
that move, merge, and are removed as the features
evolve. The performance analysis and computational-cost comparisons are
discussed in Section~\ref{sec:discussion}.

All problems are posed on $\Omega=[-L,L]^d$
with prescribed Dirichlet boundary values.
The global trial space is generated by tensor
products of the Shen--Legendre functions. The patch trial space uses the clamped Shen basis in \eqref{eq:clamped_shen_tensor_basis}.
Its corrections and their first derivatives vanish
at the patch boundary, as required by
Section~\ref{subsec:patch_compatibility}.
The fixed computational change of basis used on
the reference patch is described in
\ref{app:patch_computational_basis}.

We write $K_G$ and $K_p$ for the numbers of retained
basis functions per coordinate direction on the
global and reference patch baseplates.
Thus $r_G=K_G^d$ and $r_K=K_p^d$, with $K_p=K$
in the notation of Section~\ref{sec:methodology}.
Let $P(t)=|\Acal(t)|$ denote the number of active
patches and $P_{\max}$ its maximum during a run.
We compare four models:
\begin{enumerate}
    \item[(A)] A global solver with $K_G^d$
    coefficients.

    \item[(B)] A refined global solver with
    $K_B^d$ coefficients. We choose the smallest
    integer $K_B$ satisfying
    $K_B^d\ge K_G^d+P_{\max}K_p^d$.

    \item[(C)] Patch LegONet with $K_G^d$ global
    coefficients and $K_p^d$ coefficients per
    active patch. All experiments use $K_p=K_G$.

    \item[(D)] A physics-informed neural network
    (PINN) \cite{RaissiPINN2019}, trained using
    the governing-equation residual and the
    prescribed initial and boundary conditions,
    without solution data.
\end{enumerate}
For static patches, $P(t)$ is constant.
For moving patches, model (B) retains its fixed
coefficient count while the number of active
patches in model (C) can decrease.
The refined global baseline therefore has at least
as many coefficients as model (C) throughout
each experiment.

Models (A)--(C) use the same mechanism constructions.
Diffusion is represented by a learned diagonal
E-block fitted on the corresponding baseplate.
The Burgers experiment uses the pretrained scalar
flux $N_{\mathrm B}^{\theta}$ introduced in
Section~\ref{subsec:total_jet_assembly}.
Reaction and linear-advection operators are
supplied analytically.
Each baseplate has its own fixed numerical
operations, while physical patches with the same
reference configuration share the fitted patch
blocks. These comparisons examine global and
locally enriched representations within the same
mechanism-based framework.

Model (C) uses the initialization and state-transfer
procedure in Section~\ref{subsec:patch_management}.
The initial conditions combine localized features
with a smooth background.
The Allen--Cahn experiments use random backgrounds,
while Burgers uses a prescribed three-mode
background. The feature generators and their
parameters are specified in
\ref{app:exp_details}.
The allocation of the external source in the forced
Burgers problem is stated in
Section~\ref{subsec:exp_burgers1d}.

The reference solution is computed on a finer global
space by a spectral method with $K_{\mathrm{ref}}^d$ coefficients.
Its diffusion operator uses the computed Galerkin
eigenvalues directly, and its nonlinear functions
are evaluated analytically.
For each experiment, models (A)--(C) and the
reference solver use the same physical time step.
The integration procedure and background sampling
for the patch updates are specified in
\ref{app:patch_time_integration}.
Table~\ref{tab:exp_summary} collects the spatial
and temporal settings.

The PINN uses a common network architecture and
training procedure across the six experiments.
Following the LegONet baseline
\cite{ZhangLegONet2026}, training uses $8000$
Adam steps, with hard initial and boundary
conditions and causal weighting
\cite{WangCausalPINN2024}.
The architecture, sampling rules, and optimization
settings are given in \ref{app:pinn}.

We measure errors using the deviation
$v=u-u_\infty$, where $u_\infty$ is the constant
boundary value for the corresponding problem.
The relative error is
\begin{equation}
    e(t)
    =
    \frac{
        \|v_h(\cdot,t)-v_{\mathrm{ref}}(\cdot,t)\|_{L^2(\Omega)}
    }{
        \|v_{\mathrm{ref}}(\cdot,t)\|_{L^2(\Omega)}
    },
    \label{eq:exp_error}
\end{equation}
where $v_h=u_h-u_\infty$ is the numerical deviation
and $v_{\mathrm{ref}}=u_{\mathrm{ref}}-u_\infty$
is its reference counterpart.
The norms are evaluated on a common uniform grid with $2048$ points in 1D, $256^2$, $448^2$, and $320^2$ points for the 2D droplet, multiple-domain, and coarsening problems, and $65^3$ points in 3D.
For Burgers, $u_\infty=0$ and $v=u$.
For Allen--Cahn, subtracting the far-field state
prevents the constant background value from
dominating the normalization and obscuring errors
in the localized features.

\begin{table}[!htbp]
\centering
\caption{Discretization settings for the six experiments.
$P_{\max}$ is the maximum number of active patches.
The patch half-width $h_p$ is given in physical units,
and the capacities count basis functions per
coordinate direction.}
\label{tab:exp_summary}
\setlength{\tabcolsep}{5pt}
\begin{tabular}{llcccccccc}
\toprule
Experiment & Eq. & $L$ & $T$ & $\Delta t$
& $P_{\max}$ & $h_p$
& $K_G=K_p$ & $K_B$ & $K_{\mathrm{ref}}$ \\
\midrule
1D Allen--Cahn phase domain
& \eqref{eq:exp_ac1d}
& 6 & 5 & $5\times10^{-4}$
& 1 & 2.0 & 48 & 96 & 400 \\

1D Burgers forced shock
& \eqref{eq:exp_burgers1d}
& 6 & 5 & $5\times10^{-4}$
& 1 & 1.2 & 48 & 96 & 400 \\

2D Allen--Cahn droplet
& \eqref{eq:exp_acnd}
& 4 & 5 & $10^{-3}$
& 1 & 1.8 & 48 & 68 & 150 \\

3D Allen--Cahn droplet
& \eqref{eq:exp_acnd}
& 2.5 & 3 & $2\times10^{-3}$
& 1 & 1.4 & 28 & 36 & 60 \\

2D Allen--Cahn multiple domains
& \eqref{eq:exp_acnd}
& 7 & 5 & $10^{-3}$
& 4 & $R_i+0.55$ & 48 & 108 & 300 \\

2D advected Allen--Cahn coarsening
& \eqref{eq:exp_acdrive}
& 8 & 11 & $10^{-3}$
& 4 & adaptive & 64 & 144 & 400 \\
\bottomrule
\end{tabular}
\end{table}

All result figures use the same color convention:
blue for (A), green for (B), red for (C), orange
for the PINN (D), and black dotted lines for the
reference.
Each figure presents the initial and final fields
together with the error history $e(t)$.
Patch regions are indicated by shaded intervals
in one dimension and dashed boxes in two and
three dimensions.

\subsection{One static patch}
\label{subsec:exp_single}

We first examine local deployment with one patch
whose position and size remain fixed throughout
the run. The 1D Allen--Cahn experiment compares
global and local allocation of an equal number
of coefficients. The Burgers experiment then tests
the deployment of a pretrained nonlinear transport
mechanism. Finally, Allen--Cahn droplets in two
and three dimensions examine the effect of
dimension on the coefficient budget.

\subsubsection{Local refinement of a 1D phase domain}
\label{subsec:exp_ac1d}

The first experiment compares two ways of adding
resolution near a localized feature.
Model (B) increases the global coefficient count
from $48$ to $96$, while model (C) combines
$48$ global coefficients with $48$ patch
coefficients. We solve
\begin{equation}
    u_t=\varepsilon^2u_{xx}+u-u^3,
    \qquad
    u(\pm L,t)=1,
    \qquad
    \varepsilon=0.15,
    \label{eq:exp_ac1d}
\end{equation}
with the initial condition
\begin{equation}
    u_0(x)
    =
    1+A_0\left(
        \tanh\frac{|x|-c}{w}-1
    \right)
    +u_{\mathrm{bg}}(x),
    \qquad
    c=1.2,\quad
    w=\sqrt{2}\,\varepsilon,\quad
    A_0=0.65.
\end{equation}
The localized profile has transitions near
$x=\pm c$. The parameter $w$ sets the transition
width, and $u_{\mathrm{bg}}$ is the smooth random
background in \eqref{eq:exp_background}.
During the evolution, the interior approaches
$u\approx-1$.
A patch centered at $x=0$, with half-width
$h_p=2$, covers the localized profile and both
interfaces.

\begin{figure}[!htbp]
\centering
\includegraphics[width=\linewidth]
{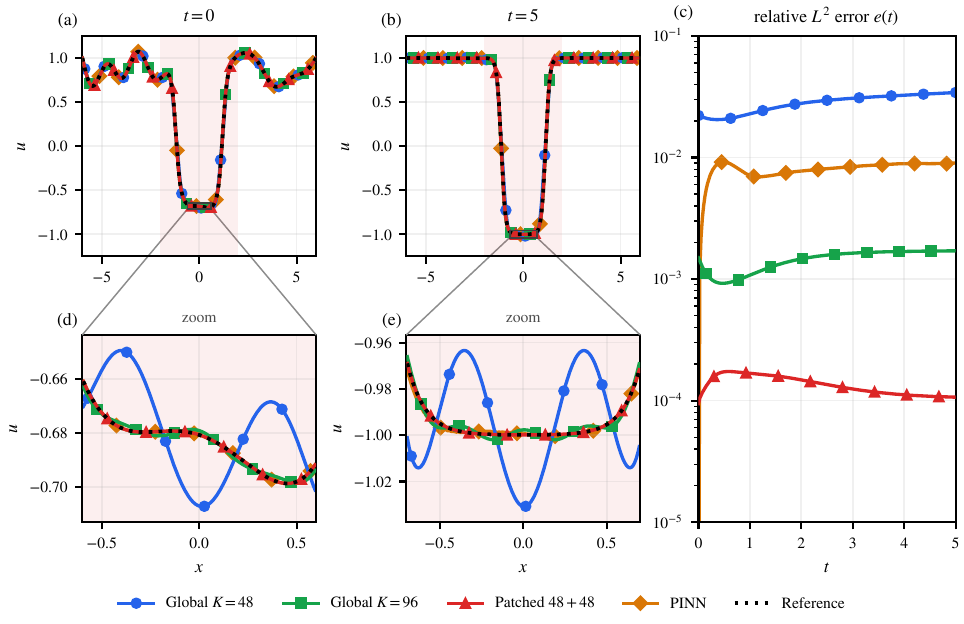}
\caption{\textbf{1D Allen--Cahn phase domain with one
static patch}, \eqref{eq:exp_ac1d}.
(a),(b)~Solutions at $t=0$ and $t=5$.
(d),(e)~Enlarged views within the phase interior.
(c)~Relative error $e(t)$ from
\eqref{eq:exp_error}.
The coefficient counts are $48$ for (A), $96$
for (B), and $48$ global plus $48$ patch
coefficients for (C); (D) is the PINN baseline with $8192$ collocation points per training step.}
\label{fig:ac1d}
\end{figure}

Models (B) and (C) have the same total coefficient
count but allocate their additional resolution
differently. Model (B) increases the global basis
dimension across the full domain.
Model (C) places the additional coefficients
around the interfaces.
The ratio $L/h_p=3$ describes the geometric
compression of its reference patch representation,
as in \eqref{eq:scale_frequency_method}.
The reference patch dimension remains fixed.

All four models reproduce the overall phase
profile in Fig.~\ref{fig:ac1d}(a),(b).
The enlarged views reveal differences within
the $u\approx-1$ phase.
Model (B) reduces the oscillations present
in the coarse global solution, while model (C)
closely follows the reference.
The oscillation amplitude in model (A) is
approximately $3\times10^{-2}$
(Fig.~\ref{fig:ac1d}(d),(e)).

The error history gives the same ordering
(Fig.~\ref{fig:ac1d}(c)).
Model (C) maintains an error between
$0.9\times10^{-4}$ and $1.7\times10^{-4}$,
with a final value of $1.1\times10^{-4}$.
This is a reduction by approximately a factor
of $16$ relative to (B) and $300$ relative
to (A).
The error of (A) increases from
$2.2\times10^{-2}$ to $3.4\times10^{-2}$.
The PINN satisfies the initial condition exactly,
but its error reaches approximately
$9\times10^{-3}$ by $t\approx0.5$ and remains
near that level.
The comparison shows the benefit of allocating
the additional coefficients near the interfaces
at a fixed total coefficient count.

\FloatBarrier
\subsubsection{Local deployment of learned Burgers transport}
\label{subsec:exp_burgers1d}

The next experiment examines local deployment
of a learned nonlinear mechanism.
We use the pretrained Burgers flux
$N_{\mathrm B}^{\theta}$ through the coupled
coefficient operator in
Section~\ref{subsec:total_jet_assembly}.
The learned parameters remain fixed, while
the patch geometry supplies the coordinate
factor from Section~\ref{subsec:analytic_pullback}. We consider the forced viscous Burgers equation
\begin{equation}
    u_t+uu_x
    =
    \nu u_{xx}+s(x),
    \qquad
    u(\pm L,t)=0,
    \qquad
    \nu=0.02.
    \label{eq:exp_burgers1d}
\end{equation}
The source is constructed from the target profile
\begin{equation}
    U^*(x)
    =
    -A\tanh\frac{x}{\delta}\,
    e^{-(x/w)^2},
    \qquad
    A=1,\quad
    \delta=0.15,\quad
    w=0.5,
\end{equation}
by setting
$s=U^*U^{*\prime}-\nu U^{*\prime\prime}$,
where primes denote differentiation with respect
to $x$.
Thus $U^*$ satisfies the stationary differential
balance.
The parameter $A$ sets the amplitude,
$\delta$ controls the central transition,
and $w$ sets the Gaussian localization width.

The initial field is
$u_0=0.3\,U^*+u_{\mathrm{bg}}$,
with the prescribed three-mode background
in \ref{app:exp_details}.
It starts from a weakened shock profile.
One patch centered at $x=0$, with half-width
$h_p=1.2$, covers the shock.
In model (C), the source is projected onto
the compatible patch space and applied only
there, as specified in
\ref{app:patch_experimental_operators}.
This treatment omits the source tail outside
the patch.

\begin{figure}[!htbp]
\centering
\includegraphics[width=\linewidth]
{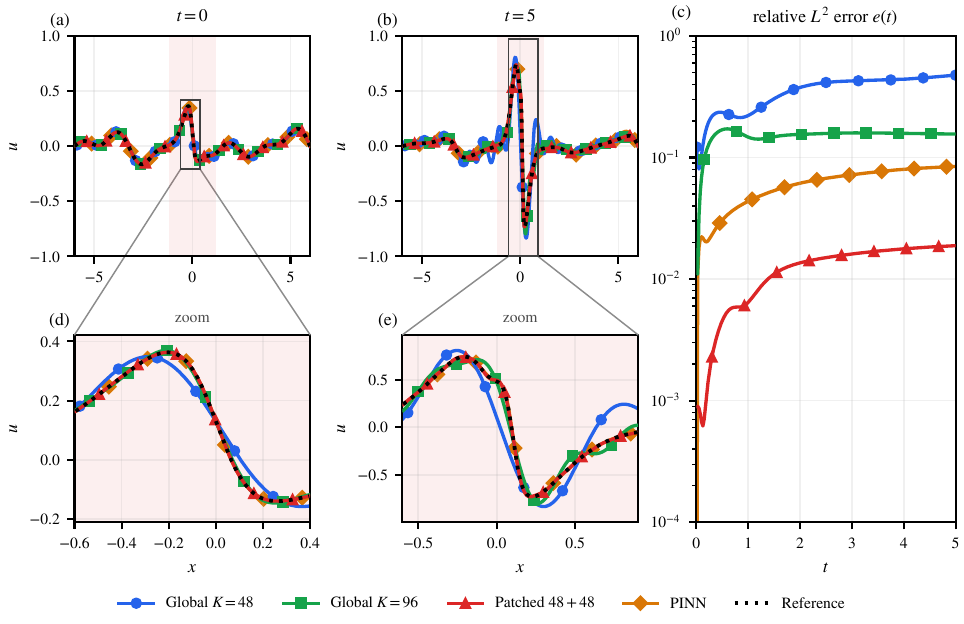}
\caption{\textbf{1D Burgers forced shock with one
static patch}, \eqref{eq:exp_burgers1d}.
(a),(b)~Solutions at $t=0$ and $t=5$.
(d),(e)~Enlarged views of the shock region.
(c)~Relative error $e(t)$ from
\eqref{eq:exp_error}.
The coefficient counts are $48$ for (A), $96$
for (B), and $48$ global plus $48$ patch
coefficients for (C); (D) is the PINN baseline with $8192$ collocation points per training step.}
\label{fig:burgers1d}
\end{figure}

The distinction between global and local
deployment becomes apparent as the shock
steepens.
During $t\lesssim1$, the error of (B) increases
from $1.7\times10^{-2}$ to approximately $0.16$
and subsequently remains near that level
(Fig.~\ref{fig:burgers1d}(c)).
Its solution exhibits oscillations on both
sides of the shock.
Model (C) resolves the transition without
the visible oscillations of the global models
(Fig.~\ref{fig:burgers1d}(b),(e)).

For the fixed patch capacity $K_p=48$,
the geometric scale ratio is $L/h_p=5$.
Model (B), in comparison, doubles the global
capacity from $48$ to $96$.
The error of (C) increases from approximately
$9\times10^{-4}$ to $1.9\times10^{-2}$
and stays below $2\times10^{-2}$ throughout
the run.
Its final error is approximately a factor
of $8$ smaller than that of (B).
Model (A) ends with an error of $0.48$,
while the PINN reaches $8.4\times10^{-2}$,
between (B) and (C).

This experiment demonstrates the use of the
same pretrained flux function at the smaller
physical scale of the patch.
The local coefficient state evolves through
the background-dependent coupling, while
the learned flux parameters remain unchanged.

\FloatBarrier
\subsubsection{Droplets in two and three dimensions}
\label{subsec:exp_ac2d3d}

The next two experiments examine one-patch
deployment in higher dimensions.
With $K_p=K_G$, model (C) uses $2K_G^d$
coefficients.
A global solver with a comparable coefficient
count can increase its capacity per direction
only to approximately $2^{1/d}K_G$.
This gives $K_B=68$ for $K_G=48$ in two
dimensions and $K_B=36$ for $K_G=28$ in
three dimensions.

We solve the Allen--Cahn equation
\begin{equation}
    u_t
    =
    \varepsilon^2\Delta u+u-u^3,
    \qquad
    u|_{\partial\Omega}=1,
    \qquad
    \varepsilon=0.1,
    \label{eq:exp_acnd}
\end{equation}
with $d=2$ and $d=3$.
The initial condition is
$u_0=1+u_{\mathrm{drop}}+u_{\mathrm{bg}}$,
where the localized droplet contains internal
oscillations and is superimposed on a smooth
random background.
The two terms are specified in
\eqref{eq:exp_droplet} and
\eqref{eq:exp_background}.
The droplet is centered at the origin,
with radius $R=1.26$ in two dimensions
and $R=0.9$ in three dimensions.
A single patch covers it, with half-width
$1.8$ and $1.4$, respectively.

\begin{figure}[!htbp]
\centering
\includegraphics[width=\linewidth]
{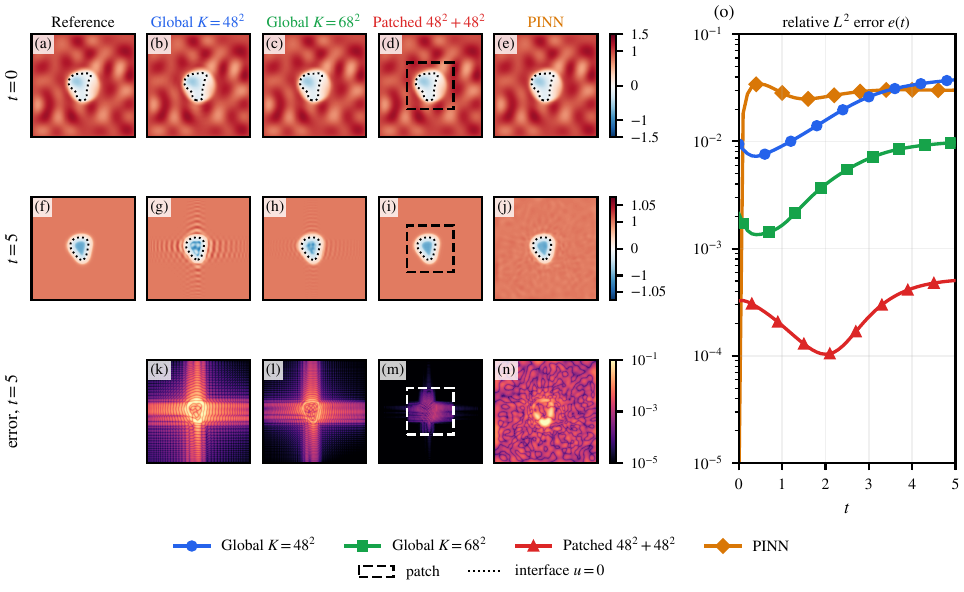}
\caption{\textbf{2D Allen--Cahn droplet with one
static patch}, \eqref{eq:exp_acnd}.
Columns show the reference, (A) $48^2$
coefficients, (B) $68^2$ coefficients,
(C) $48^2$ global plus $48^2$ patch
coefficients, and (D) the PINN with $16384$ collocation points per training step.
(a)--(e)~Solutions at $t=0$.
(f)--(j)~Solutions at $t=5$; dotted lines
mark the interface $u=0$.
(k)--(n)~Pointwise error
$|u_h-u_{\mathrm{ref}}|$ at $t=5$
on a logarithmic scale.
(o)~Relative error $e(t)$ from
\eqref{eq:exp_error}.}
\label{fig:ac2d}
\end{figure}

\begin{figure}[!htbp]
\centering
\includegraphics[width=\linewidth]
{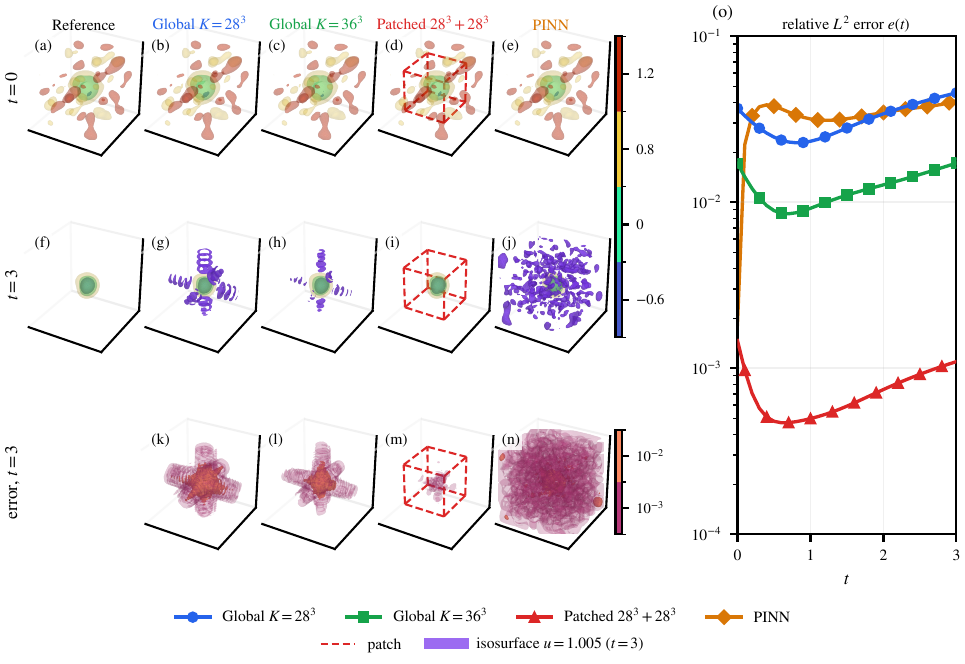}
\caption{\textbf{3D Allen--Cahn droplet with one
static patch}, \eqref{eq:exp_acnd}.
Columns show the reference, (A) $28^3$
coefficients, (B) $36^3$ coefficients,
(C) $28^3$ global plus $28^3$ patch
coefficients, and (D) the PINN with $32768$ collocation points per training step.
(a)--(e)~Isosurfaces at $t=0$.
(f)--(j)~Isosurfaces at $t=3$; the violet
surface corresponds to $u=1.005$,
slightly above the boundary value, and reveals background oscillations.
(k)--(n)~Isosurfaces of
$|u_h-u_{\mathrm{ref}}|$ at $t=3$,
at levels $10^{-3}$ (transparent)
and $10^{-2}$ (opaque).
(o)~Relative error $e(t)$ from
\eqref{eq:exp_error}.}
\label{fig:ac3d}
\end{figure}

In two dimensions, the droplet interior approaches
the $u\approx-1$ phase while its interface sharpens
(Fig.~\ref{fig:ac2d}).
The error of (B) increases from
$2.1\times10^{-3}$ to $9.8\times10^{-3}$.
Model (C) maintains an error between
$1.0\times10^{-4}$ and $5.0\times10^{-4}$,
with a final error approximately a factor
of $19$ smaller than that of (B).
The coarse global model (A) increases from
$9.5\times10^{-3}$ to $3.8\times10^{-2}$.

The pointwise error plots also distinguish
the representations
(Fig.~\ref{fig:ac2d}(k)--(m)).
The global models exhibit cross-shaped error
patterns along the coordinate directions.
For model (C), the visible errors are more
localized around the droplet.
The patch basis functions are supported
on $[-h_p,h_p]^2$, while the global component
continues to represent the surrounding
smooth field.

In three dimensions, the final errors are
$1.7\times10^{-2}$ for (B) and
$1.1\times10^{-3}$ for (C), giving a reduction
by approximately a factor of $16$.
The error of (A) is $4.6\times10^{-2}$
(Fig.~\ref{fig:ac3d}).
The comparable coefficient budget increases
the global capacity per direction by
$36/28\approx1.3$.
The patch instead has the geometric scale
ratio $L/h_p\approx1.8$ at its fixed
reference capacity.

The PINN errors are $3.0\times10^{-2}$ in
two dimensions and $4.1\times10^{-2}$ in
three dimensions, comparable to model (A).
Its pointwise errors extend across both
the background and droplet regions
(Figs.~\ref{fig:ac2d}(n)
and~\ref{fig:ac3d}(n)).
In both dimensions, the patched solver
achieves the lowest final error despite
using slightly fewer coefficients than
the refined global baseline.
\FloatBarrier

\subsection{Multiple static patches}
\label{subsec:exp_multi}

This experiment examines the simultaneous deployment
of the same reference patch blocks around several
localized features. Each patch carries its own
coefficient state, while a single global background
represents the surrounding field.

We solve the two-dimensional Allen--Cahn equation
\eqref{eq:exp_acnd} with four separated droplets
on a smooth random background:
$u_0=1+\sum_{i=1}^{4}u_{\mathrm{drop},i}
+u_{\mathrm{bg}}$.
The droplet profiles and background are defined
in \eqref{eq:exp_droplet}
and \eqref{eq:exp_background}.
The droplet radii are approximately one, with
their centers and radii listed in
Table~\ref{tab:exp_multi}.
Each droplet is covered by a patch centered
at the same location, with half-width
$R_i+0.55$.
The four patches have nonoverlapping interiors
and remain fixed throughout the run.

Model (C) uses $48^2$ global coefficients
and four patch states with $48^2$ coefficients
each, giving $5\cdot48^2=11520$ coefficients
in total.
All four patches share the same frozen
reference patch library.
The refined global model (B) uses
$108^2=11664$ coefficients, slightly more
than the patched solver.

\begin{figure}[!htbp]
\centering
\includegraphics[width=\linewidth]
{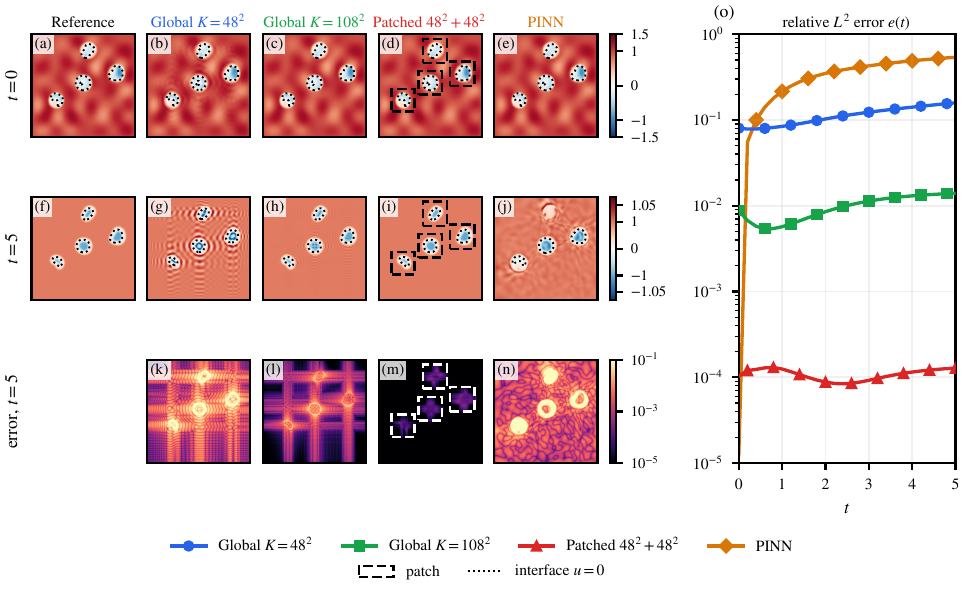}
\caption{\textbf{2D Allen--Cahn multiple domains
with four static patches},
\eqref{eq:exp_acnd}.
Columns show the reference,
(A) $48^2$ coefficients,
(B) $108^2$ coefficients,
(C) $48^2$ global coefficients plus four
patch states with $48^2$ coefficients each,
and (D) the PINN with $16384$ collocation points per training step.
(a)--(e)~Solutions at $t=0$.
(f)--(j)~Solutions at $t=5$; dotted lines
mark the interfaces $u=0$.
(k)--(n)~Pointwise error
$|u_h-u_{\mathrm{ref}}|$ at $t=5$
on a logarithmic scale.
(o)~Relative error $e(t)$ from
\eqref{eq:exp_error}.}
\label{fig:ac2d4}
\end{figure}

Model (C) maintains an error between
$0.8\times10^{-4}$ and $1.3\times10^{-4}$
throughout the run
(Fig.~\ref{fig:ac2d4}(o)).
Its error remains of the same order of
magnitude as in the single-droplet experiment.
The error of (B) increases from
$8.8\times10^{-3}$ to $1.4\times10^{-2}$.
At $t=5$, the error reduction relative
to (B) is approximately a factor of $109$,
compared with approximately $19$ in the
single-droplet test.

The pointwise error plots show different
spatial patterns.
The global models exhibit cross-shaped
errors around the four droplets, extending
along the coordinate directions
(Fig.~\ref{fig:ac2d4}(k),(l)).
For model (C), the visible errors are more
localized around the individual features
(Fig.~\ref{fig:ac2d4}(m)).
The patch sizes give geometric compression
ratios $L/h_p\approx4.4$, while global
refinement increases the capacity per
direction by $108/48=2.25$.
The additional patch coefficients are
therefore concentrated near the droplets,
with the global state continuing to
represent the smooth background.

The PINN has a final error of $0.54$,
compared with $3.0\times10^{-2}$ in the
single-droplet experiment.
Its solution omits one droplet and
misplaces the interfaces of the others
(Fig.~\ref{fig:ac2d4}(j),(n)).
The patched solution retains all four
droplets and closely follows their
reference interfaces at the comparable
coefficient budget.
\FloatBarrier

\subsection{Moving patches with merging and removal}
\label{subsec:exp_moving}

The final experiment examines evolving features
that require patch relocation, merging, and removal.
The coefficient states must be transferred between
successive configurations while the pretrained
reference blocks remain fixed. We consider two-dimensional advected Allen--Cahn coarsening,
\begin{equation}
    u_t+\bm V\cdot\nabla u
    =
    \varepsilon^2\Delta u+u-u^3+H,
    \qquad
    u|_{\partial\Omega}=u_+,
    \label{eq:exp_acdrive}
\end{equation}
with $\varepsilon=0.1$, $H=-0.33$, and
$\bm V=(0.1,0.05)$.
The constant velocity $\bm V$ transports the
features, while $H$ favors the negative phase.
The boundary value $u_+\approx0.7472$ is the
largest root of $u-u^3+H=0$.
The initial condition contains three droplets
of radius $0.6$ and one smaller droplet
of radius $0.18$ on a smooth random background.
Their locations and profile parameters are
listed in Table~\ref{tab:exp_moving}.
The larger droplets grow, two eventually merge,
and the smaller droplet shrinks and disappears.
Linear advection is specified analytically
on each baseplate, with the patch factor
$h_p^{-1}$ from
Section~\ref{subsec:analytic_pullback}.

The patch configuration is selected from the
current numerical field using the procedure
in Section~\ref{subsec:patch_management}.
Every $100$ time steps, the zeroth-order
feature monitor \eqref{eq:truncation_monitor}
is evaluated with $M=0$.
The one-sided criterion $u_h<u_+-0.45$
\eqref{eq:patch_phase_detection} identifies the droplet phase.
Square patches cover the detected regions,
and nearby candidates are merged.
A patch is retained after its feature disappears
until its correction satisfies the removal
tolerance.
The detection and geometric parameters are
specified in \ref{app:exp_moving}.

Patch centers and half-widths remain fixed
between rebuilding events.
When the configuration changes, the current
field is redistributed using
\eqref{eq:patch_rebuild_local}--\eqref{eq:patch_rebuild_global}.
The new states are computed from the old total
field and background, with transfer accuracy
limited by the corresponding projections.
The reference bases and learned mechanism
parameters remain unchanged.

\begin{figure}[!htbp]
\centering
\includegraphics[width=\linewidth]
{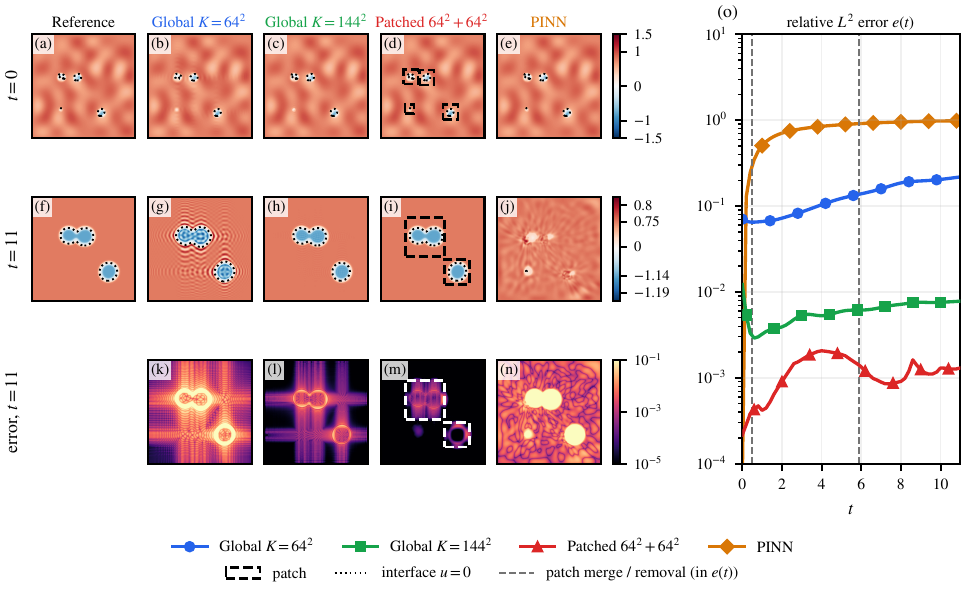}
\caption{\textbf{2D advected Allen--Cahn coarsening
with moving patches},
\eqref{eq:exp_acdrive}.
Columns show the reference,
(A) $64^2$ coefficients,
(B) $144^2$ coefficients,
(C) $64^2$ global coefficients plus
$64^2$ coefficients per active patch,
and (D) the PINN with $16384$ collocation points per training step.
(a)--(e)~Solutions at $t=0$;
model (C) has four patches.
(f)--(j)~Solutions at $t=11$;
model (C) has two patches.
Dotted lines mark the interfaces $u=0$,
and dashed boxes indicate the patch regions
in model (C).
(k)--(n)~Pointwise error
$|u_h-u_{\mathrm{ref}}|$ at $t=11$
on a logarithmic scale.
(o)~Relative error $e(t)$ from
\eqref{eq:exp_error}.
Vertical dashed lines mark the merger
of the patches associated with droplets
1 and 2 at $t=0.5$ and the removal
of the patch associated with droplet 4
at $t=5.9$.}
\label{fig:acdrive}
\end{figure}

The configuration is rebuilt ten times during
the run.
The patches associated with droplets 1 and 2
merge at $t=0.5$, and the patch associated
with droplet 4 is removed at $t=5.9$.
Eight further updates adjust patch positions
or sizes as the features drift and grow.
The active patch count decreases from four
to two.

The error of (C) remains below
$2.1\times10^{-3}$ throughout the run
and ends at $1.3\times10^{-3}$
(Fig.~\ref{fig:acdrive}(o)).
No sustained increase is visible following
the rebuilding events.
Model (B) has a final error of
$7.8\times10^{-3}$, giving an error reduction
by approximately a factor of $6$ for
the patched solver.

This accuracy is obtained with a decreasing
coefficient count.
Model (B) retains $144^2=20736$ coefficients
throughout the run.
Model (C) begins with
$5\cdot64^2=20480$ coefficients and ends
with $3\cdot64^2=12288$.
Its final state therefore uses approximately
$59\%$ of the coefficients in (B).

The error of the coarse global model (A)
increases to $0.22$.
The PINN reaches an error of $0.99$ at
$t=11$, with the droplet structures largely
absent from its predicted field
(Fig.~\ref{fig:acdrive}(j)).
The patched solver continues to represent
the evolving features through changes
in geometry and coefficient states,
using the same frozen reference patch blocks.

\FloatBarrier

\FloatBarrier
\subsection*{\textbf{Summary of numerical results}}

Table~\ref{tab:exp_results} summarizes the final
relative errors across the six experiments in
Sections~\ref{subsec:exp_single}--\ref{subsec:exp_moving}.
Patch LegONet (C) achieves the lowest final error
in each case.
The refined global baseline (B) uses at least
as many coefficients as the patched solver
throughout each run.
Section~\ref{sec:discussion} discusses these
results and examines the computational cost.

\begin{table}[!htbp]
\centering
\caption{Final relative $L^2$ errors
$e(T)$ from \eqref{eq:exp_error}
for the six experiments.
Baseline (B) uses at least as many coefficients
as Patch LegONet (C).
The last column reports the approximate
error-reduction factor relative to (B).}
\label{tab:exp_results}
\setlength{\tabcolsep}{5pt}
\begin{tabular}{lccccc}
\toprule
Experiment
& (A) & (B) & (C) & (D) PINN & (B)/(C) \\
\midrule
1D Allen--Cahn phase domain
& $3.4\times10^{-2}$
& $1.7\times10^{-3}$
& $\bm{1.1\times10^{-4}}$
& $9.0\times10^{-3}$
& 16 \\

1D Burgers forced shock
& $4.8\times10^{-1}$
& $1.6\times10^{-1}$
& $\bm{1.9\times10^{-2}}$
& $8.4\times10^{-2}$
& 8 \\

2D Allen--Cahn droplet
& $3.8\times10^{-2}$
& $9.8\times10^{-3}$
& $\bm{5.0\times10^{-4}}$
& $3.0\times10^{-2}$
& 19 \\

3D Allen--Cahn droplet
& $4.6\times10^{-2}$
& $1.7\times10^{-2}$
& $\bm{1.1\times10^{-3}}$
& $4.1\times10^{-2}$
& 16 \\

2D Allen--Cahn multiple domains
& $1.6\times10^{-1}$
& $1.4\times10^{-2}$
& $\bm{1.3\times10^{-4}}$
& $5.4\times10^{-1}$
& 109 \\

2D advected Allen--Cahn coarsening
& $2.2\times10^{-1}$
& $7.8\times10^{-3}$
& $\bm{1.3\times10^{-3}}$
& $9.9\times10^{-1}$
& 6 \\
\bottomrule
\end{tabular}
\end{table}

\FloatBarrier

\section{Discussion}
\label{sec:discussion}

The experiments show how a fixed reference patch
library can provide additional physical resolution
through its deployment. The main findings concern
reuse across scales, accuracy at comparable
coefficient counts, and the cost of following
localized features during evolution.

\textbf{Reuse across physical scales.}
Patch LegONet separates the reference representation
of a mechanism from its physical deployment.
For a given reference patch basis, capacity, and
dimension, physical patches share the same frozen
blocks while carrying independent coefficient states.
The diffusion block is fitted once for that reference
configuration.
For Burgers transport, the pretrained scalar flux
is retained, with fixed spatial operations constructed
for the chosen baseplate.
The global baseplate uses blocks compatible with
its own representation.
This construction extends mechanism reuse to
different physical locations and scales without
enlarging the reference patch blocks.

\textbf{Accuracy at comparable coefficient counts.}
Patch LegONet achieves the lowest final error
in all six experiments, including comparisons
with a globally refined solver that uses at least
as many coefficients
(Table~\ref{tab:exp_results}).
The error reduction relative to (B) is approximately
one order of magnitude in most tests.
The pointwise error plots also show more localized
errors around the patched features, compared with
the extended patterns visible in the global solutions.
These results support allocating additional
coefficients near localized structures while
retaining a global representation of the
surrounding field.

\textbf{Multiple features and higher dimensions.}
The improvement persists in the higher-dimensional
experiments, where the coefficient count of a
tensor-product global representation grows as $K^d$.
For a single Allen--Cahn feature, the error reduction
relative to (B) is approximately a factor of
$16$--$19$ from one to three dimensions.
The four-droplet experiment gives the largest
observed reduction, approximately a factor of $109$.
In that configuration, separate patch states
represent the four features using one shared
reference patch library.
The results illustrate local deployment for
several simultaneous structures as well as
for features in higher dimensions.

\textbf{Patch updates without retraining.}
The moving experiment changes the patch configuration
ten times through relocation, resizing, merging,
and removal.
The error remains below $2.1\times10^{-3}$
throughout the run and ends approximately a
factor of $6$ below that of (B).
After merging and removal, the patched state
contains approximately $59\%$ of the coefficients
in the refined global baseline.
The geometry and coefficient states change
through the transfer procedure in
Section~\ref{subsec:patch_management},
while the reference patch blocks remain fixed.
The observed error history supports this use
of state transfer for the evolving features
considered here.

\phantomsection
\label{subsec:exp_summary}
\textbf{Computational cost.}
At matched final error, Patch LegONet reduces
runtime by approximately a factor of $1.7$
for 1D Burgers and $2.8$--$7.0$ for the
two-dimensional experiments
(Table~\ref{tab:exp_cost}).
The global resolution is selected to match
the final error of (C) within $15\%$, and
both coefficient-based solvers are timed
on one CPU thread.
These global capacities can differ from
those of baseline (B) in Section~\ref{sec:experiments}.
For 1D Allen--Cahn, the runtimes are nearly
equal at $0.82$ and $0.80$~s.
In three dimensions, the largest tested global
configuration has $52^3$ coefficients and
takes more than $355$~s, compared with
$106$~s for (C), while its final error
remains $2.3$ times larger.
The three-dimensional result is therefore
reported as a separate comparison without
matched final error.

\begin{table}[!htbp]
\centering
\caption{Wall times in seconds for global
refinement and Patch LegONet on one CPU thread.
The global resolution is selected to match
the final error of (C) within $15\%$,
except for the marked three-dimensional case.
PINN training times on one A100 GPU are
reported separately.
$^\dagger$The largest tested 3D global
configuration has $52^3$ coefficients;
its final error remains $2.3$ times
that of (C).}
\label{tab:exp_cost}
\setlength{\tabcolsep}{8pt}
\begin{tabular}{lccc}
\toprule
Experiment
& Global refinement
& Patch LegONet
& PINN training \\
\midrule
1D Allen--Cahn phase domain
& 0.82
& 0.80
& 197 \\

1D Burgers forced shock
& 97
& \textbf{57}
& 219 \\

2D Allen--Cahn droplet
& 20
& \textbf{7.2}
& 478 \\

3D Allen--Cahn droplet
& $>355^\dagger$
& \textbf{106}
& 915 \\

2D Allen--Cahn multiple domains
& 154
& \textbf{22}
& 333 \\

2D advected Allen--Cahn coarsening
& 337
& \textbf{95}
& 397 \\
\bottomrule
\end{tabular}
\end{table}

\textbf{Comparison with the PINN baseline.}
The specified PINN configuration has substantially
larger errors in the multiple-feature and
coarsening experiments, reaching $0.54$
and $0.99$, respectively.
The patched solver continues to represent
the localized structures accurately in these tests.
This comparison uses the architecture
and training procedure in \ref{app:pinn}.
\FloatBarrier

\section{Conclusion}
\label{sec:conclusion}

We introduced Patch LegONet to bring local refinement
into PDE solvers assembled from reusable mechanism
blocks. Deploying low-frequency reference blocks
on smaller physical patches provides higher local
resolution while keeping the pretrained parameters
fixed. The numerical experiments demonstrate better
accuracy than global refinement at comparable
coefficient counts, together with runtime savings
in several matched-error comparisons. This reuse
also supports changing patch configurations as
localized features evolve.

The broader implication is that a pretrained
mechanism library need not remain tied to a single
physical resolution. Spatial deployment provides
an additional way to allocate approximation
capacity according to the solution, extending
the usefulness of existing reference blocks
to localized multiscale dynamics. Future work
can explore hierarchical patch deployment
and extensions to more general geometries.

\section*{Declaration of generative AI and AI-assisted
technologies in the manuscript preparation process}

During the preparation of this manuscript, the authors
used ChatGPT (OpenAI) to assist with language revision,
manuscript organization, and identification of relevant
literature. The authors reviewed and edited the resulting
content, verified the references, and take full
responsibility for the final manuscript.

\section*{Acknowledgment}
We gratefully acknowledge the support from DARPA MATHBAC program under grant number: HR0011-26-3-E152, the National Science Foundation (DMS-2533878, DMS-2053746, DMS-2134209, ECCS-2328241, CBET-2347401 and OAC-2311848), and U.S.~Department of Energy (DOE) Office of Science Advanced Scientific Computing Research program under the "Uncertainty Quantification for Multifidelity Operator Learning (MOLUcQ)" project (Project No. 81739), DE-SC0023161, the SciDAC LEADS Institute, and DOE–Fusion Energy Science, under grant number: DE-SC0024583.

\bibliographystyle{elsarticle-num}
\bibliography{mybib}

\newpage
\appendix

\section{Reference block implementation and time integration}
\label{app:patch_implementation}

The method in Section~\ref{sec:methodology} is formulated
for a general boundary-compatible basis. This appendix
specifies the fixed computational coordinates, diffusion
block, and time integrator used in the experiments.

{
\subsection{Computational basis}
\label{app:patch_computational_basis}

Let $B_0$ and $D_0$ collect the values and derivatives
of the one-dimensional clamped Shen basis in
\eqref{eq:clamped_shen_basis} at the reference quadrature
nodes. With quadrature weights $W$, form
\begin{equation}
    M_0=B_0^{\mathsf T}WB_0,
    \qquad S_0=D_0^{\mathsf T}WD_0.
    \label{eq:patch_raw_matrices}
\end{equation}
We solve the generalized eigenproblem
\begin{equation}
    S_0 V=M_0 V\Lambda_K,
    \qquad V^{\mathsf T}M_0V=I,
    \qquad \Lambda_K=\operatorname{diag}(\lambda_k).
    \label{eq:patch_reference_eigenbasis}
\end{equation}
The computational basis has nodal and derivative matrices
$B=B_0V$ and $D=D_0V$, giving
$M_K=B^{\mathsf T}WB=I$ and
$D^{\mathsf T}WD=\Lambda_K$.
It spans the same trial space and retains all homogeneous
boundary conditions. If $a_0$ denotes coefficients in the
original basis, then $a_0=Va$ represents the same field
in the computational coordinates $a$. In multiple
dimensions, we use tensor products of these transformed
one-dimensional bases; the reference Laplacian eigenvalues
are the corresponding sums of the one-dimensional values.

This change of basis is performed once when constructing
the reference baseplate, before training the diffusion
block. Reconstruction, projection, and mechanism
evaluation then use these same coordinates throughout
training and deployment. In the implementation formulas
below, $a_p$, $B$, and $D$ refer to these computational
coordinates, and the projection is $\Phat=B^{\mathsf T}W$.
The general formulas in Section~\ref{sec:methodology}
remain valid with this choice. Quadrature uses $3K+12$ (patch) and $3K+8$ (global) Gauss--Legendre nodes per direction.}

\subsection{Diagonal diffusion block}
\label{app:patch_diffusion}

{In the computational basis defined above, we fit a
quadratic E-block following LegONet's energy-generator
construction:
\begin{equation}
    E_K^\theta(a)=\frac12 a^{\mathsf T}\Lambda_K^\theta a,
    \qquad
    \Lambda_K^\theta=\operatorname{diag}(e^{\theta_k}),
    \qquad
    F_{\mathrm{diff},K}^\theta(a)
    =-\nabla_a E_K^\theta(a)=-\Lambda_K^\theta a.
    \label{eq:patch_diffusion_construction}
\end{equation}
The training targets are the instantaneous reference
rates $-\Lambda_Ka$ for sampled
coefficient states, where $\lambda_k$ are the reference
Galerkin eigenvalues. This fits one decay rate per
coefficient. On physical patch $p$, the diffusion
contribution is simply
$-\nu h_p^{-2}\Lambda_K^\theta a_p$.
Neither diffusion training nor the eigenproblem is
repeated for individual physical patches. This diagonal
quadratic E-block differs from the original one-dimensional
LegONet Burgers implementation, which extracts a dense
stiffness matrix from its trained energy network.
The reference patch fit is specific to its basis,
capacity, and dimension, and is shared across physical
patches with that reference configuration. All the
experiments in Section~\ref{sec:experiments} use this
diffusion construction.}

\subsection{Nodal realization of reaction and transport}
\label{app:patch_nodal_realization}

{
For reference quadrature weights
$W=\operatorname{diag}(w_q)$, the coefficient projection is
$\Phat=M_K^{-1}B^{\mathsf T}W$, with
$M_K=B^{\mathsf T}WB$. This formula applies to a general
compatible basis; the computational basis in~\ref{app:patch_computational_basis} gives $M_K=I$.
The reaction contribution is therefore
$M_K^{-1}B^{\mathsf T}W[f(g+v)-f(g)]$.
The Allen--Cahn experiments use the closed-form reaction
$f(u)=u-u^3$, with the constant boundary lifting included
in $g$. Equivalently, the code evolves the field relative
to the lifting and evaluates the shifted reaction.
For the driven reaction $f(u)=u-u^3+H$, the constant $H$
cancels in the patch difference but remains in the
global reaction.

The original LegONet Burgers implementation
\cite{ZhangLegONet2026} realizes transport through
\begin{equation}
    a\longmapsto B_{\mathrm L}a
     \longmapsto N_{\mathrm B}^{\theta}(B_{\mathrm L}a)
     \longmapsto P_{\mathrm L}N_{\mathrm B}^{\theta}(B_{\mathrm L}a)
     \longmapsto J_{\mathrm L}P_{\mathrm L}
                          N_{\mathrm B}^{\theta}(B_{\mathrm L}a).
    \label{eq:legonet_transport_path}
\end{equation}
Here $B_{\mathrm L}$ reconstructs the field,
$P_{\mathrm L}$ projects onto the original coefficient
space, and $J_{\mathrm L}$ represents the reference
derivative; the PDE assembly supplies the transport sign.
We retain the learned density derivative
$N_{\mathrm B}^{\theta}$ and realize the derivative
on the compatible patch basis. Define
$D_{qk}=\widehat\phi_k'(\xi_q)$ and
$\delta N_q=N_{\mathrm B}^{\theta}(g_q+v_q)
             -N_{\mathrm B}^{\theta}(g_q)$.
Integration by parts in \eqref{eq:patch_block} gives
\begin{equation}
    \Ccal_{\mathrm{adv},p}
    =\frac{1}{h_p}M_K^{-1}D^{\mathsf T}W\delta N,
    \qquad
    \widehat{\mathcal Q}_{\mathrm{adv},K}^{\theta}(z)
    =M_K^{-1}D^{\mathsf T}W N_{\mathrm B}^{\theta}(z).
    \label{eq:patch_transport_weak_form}
\end{equation}
Here $\widehat{\mathcal Q}_{\mathrm{adv},K}^{\theta}$ includes the transport sign. The boundary term vanishes because the patch test basis
is zero at both endpoints. For the exact quadratic flux,
$\delta N=gv+v^2/2$, with componentwise products.
The original coefficient derivative matrix belongs to
its original baseplate and is not copied to the patch
basis. The reference patch matrices are constructed
once and shared by all physical patches with that
reference configuration.}

\subsection{Assembled experimental operators}
\label{app:patch_experimental_operators}

{In the forced Burgers experiment, the source
is assigned only to the patch, with $\widehat s_p=s\circ T_p$;
the global background receives no copy of it.
The diffusion and transport blocks together give the
patch equation used in the Burgers experiment of Section~\ref{subsec:exp_burgers1d}:
\begin{equation}
    \dot a_p
    =-\frac{\nu}{h_p^2}\Lambda_K^\theta a_p
     +\frac{1}{h_p}M_K^{-1}D^{\mathsf T}W
        \left[N_{\mathrm B}^{\theta}(g+Ba_p)
              -N_{\mathrm B}^{\theta}(g)\right]
     +\Phat[\widehat s_p].
    \label{eq:burgers_complete_patch_block}
\end{equation}
The frozen learned quantities are the reference decay
rates in $\Lambda_K^\theta$ and the weights defining
$N_{\mathrm B}^{\theta}$. The matrices $B,D,W,M_K$
are fixed reference numerical data. Each deployed
patch supplies its own $a_p$, $h_p$, and sampled
background $g$; a numerical integrator advances $a_p$
using this right-hand side. Together, the frozen mechanisms, reference matrices,
and patch-specific inputs define the patch block. For the Allen--Cahn experiments
in Sections~\ref{subsec:exp_ac1d}, \ref{subsec:exp_ac2d3d},
\ref{subsec:exp_multi}, and \ref{subsec:exp_moving},
the reaction contribution instead uses
$\Phat[f(g+Ba_p)-f(g)]$ from
\eqref{eq:patch_reaction_block}, with diffusion treated
in the same way. The gain in local spatial resolution comes from
evaluating the fixed reference basis at
$T_p^{-1}(\x)$, while the equations above specify how
its coefficients evolve.}

\subsection{Time integration}
\label{app:patch_time_integration}

{The time integrator used in Section~\ref{sec:experiments}
composes a diffusion half-step, a full reaction/transport
step, and another diffusion half-step
\cite{Strang1968,Hairer2006}. Let $t_n=n\Dt$.
We first advance $a_G^n$ to $a_G^{n+1}$ and sample
$g_p^n=u_G^n\circ T_p$ and
$g_p^{n+1}=u_G^{n+1}\circ T_p$ at the patch nodes.
For a patch, the diffusion block multiplies each
coefficient by $\exp[-\nu\lambda_k^\theta\Dt/(2h_p^2)]$
in each half-step, where
$\lambda_k^\theta=(\Lambda_K^\theta)_{kk}$ is its
frozen learned decay rate from
\eqref{eq:patch_diffusion_construction}. Between these half-steps,
Heun's method evaluates the reaction/transport rate
first at the diffused patch state with background
$g_p^n$, then at the Euler-predicted patch state with
background $g_p^{n+1}$, and averages the two rates.
The local update thus accounts for the changing global
background while the patch geometry remains fixed. The forced Burgers source is added
after this composition by a first-order Euler step. For the stationary source in that
experiment, its coefficient vector is precomputed and added
as $\Dt\,\Phat[\widehat s_p]$; the global background
receives no copy of this source.}

\section{Details of the numerical experiments}
\label{app:exp_details}

Allen--Cahn problems are solved for $v=u-u_\infty$, so that
$v|_{\partial\Omega}=0$. All random draws use fixed seeds.

\subsection{Initial conditions}
The background is
\begin{equation}
    u_{\mathrm{bg}}(\x)=s\sum_{i=1}^{M}c_i\prod_{j=1}^{d}
    \sin\frac{n_{ij}\pi(x_j+L)}{2L},
    \qquad n_{ij}\sim\mathcal U\{1,\dots,N\},\quad
    c_i\sim\frac{\mathcal N(0,1)}{\sqrt{\sum_j n_{ij}}},
    \label{eq:exp_background}
\end{equation}
with $s$ chosen so that $\max|u_{\mathrm{bg}}|=A_{\mathrm{bg}}$. The 1D Burgers
problem instead uses the fixed background
$u_{\mathrm{bg}}(x)=\sum_{(c,n)\in\mathcal S}c\,\sin\frac{n\pi(x+L)}{2L}$,
$\mathcal S=\{(0.08,5),(-0.06,8),(0.05,13)\}$.

A droplet centered at $\bm c$ with radius $R$ is
\begin{equation}
    u_{\mathrm{drop}}(\x)=m(\x)\big(g(\x-\bm c)-u_\infty\big),\qquad
    m(\x)=\tfrac12\Big(1-\tanh\frac{|\x-\bm c|-R}{w}\Big),
    \label{eq:exp_droplet}
\end{equation}
where $g(\bm y)=\sum_{j=1}^{J}\cos(k_j\bm\theta_j\cdot\bm y+\varphi_j)$,
$|\cdot|$ is the Euclidean norm, $k_j\sim\mathcal U[1.2/\sigma,2.6/\sigma]$, $\sigma=R/1.8$, $\bm\theta_j$ uniform
on the unit sphere, $\varphi_j\sim\mathcal U[0,2\pi)$, rescaled affinely to
$[g_{\min},g_{\max}]$ on $|\bm y|\le R$.

\subsection{Parameters}
Table~\ref{tab:exp_acnd} lists the initial-condition parameters of the
Allen--Cahn experiments; the droplets of the 2D Allen--Cahn multiple domains and the 2D
advected Allen--Cahn coarsening are listed in
Tables~\ref{tab:exp_multi} and~\ref{tab:exp_moving}.

\begin{table}[!htbp]
\centering
\caption{Initial-condition parameters of the Allen--Cahn problems; $w_i=R_i/7.2$ per droplet.}
\label{tab:exp_acnd}
\begin{tabular}{lccccc}
\toprule
Experiment & $R$ & $w$ & $J$ & $[g_{\min},g_{\max}]$ & $(M,N,A_{\mathrm{bg}})$ \\
\midrule
1D phase domain   & --   & --    & -- & --           & $(16,16,0.4)$ \\
2D droplet        & 1.26 & 0.175 & 6  & $[-0.3,0.5]$ & $(24,10,0.4)$ \\
3D droplet        & 0.9  & 0.125 & 8  & $[-0.6,0.2]$ & $(30,6,0.4)$  \\
multiple domains  & Table~\ref{tab:exp_multi} & $R_i/7.2$ & 6 & $[-0.6,0.2]$ & $(24,10,0.4)$ \\
coarsening        & Table~\ref{tab:exp_moving} & 0.15 & 6 & Table~\ref{tab:exp_moving} & $(24,10,0.3)$ \\
\bottomrule
\end{tabular}
\end{table}

\begin{table}[!htbp]
\centering
\caption{Droplets of the 2D Allen--Cahn multiple domains experiment. Patch $i$ is centered at
$\bm c_i$ with half-width $R_i+0.55$.}
\label{tab:exp_multi}
\begin{tabular}{ccc}
\toprule
$i$ & $\bm c_i$ & $R_i$ \\
\midrule
1 & $(0.02,\ 0.39)$   & 1.049 \\
2 & $(-3.69,\ -1.96)$ & 0.985 \\
3 & $(4.35,\ 1.64)$   & 1.113 \\
4 & $(0.64,\ 4.72)$   & 1.089 \\
\bottomrule
\end{tabular}
\end{table}

\begin{table}[!htbp]
\centering
\caption{Initial droplets of the 2D advected Allen--Cahn coarsening experiment
\eqref{eq:exp_acdrive}.}
\label{tab:exp_moving}
\begin{tabular}{cccc}
\toprule
$i$ & $\bm c_i$ & $R_i$ & $[g_{\min},g_{\max}]$ \\
\midrule
1 & $(-3.4,\ 1.4)$  & 0.6  & $[-0.6,0.2]$ \\
2 & $(-0.9,\ 1.4)$  & 0.6  & $[-0.6,0.2]$ \\
3 & $(2.7,\ -4.0)$  & 0.6  & $[-0.6,0.2]$ \\
4 & $(-3.5,\ -3.5)$ & 0.18 & $[-0.3,0.5]$ \\
\bottomrule
\end{tabular}
\end{table}

\subsection{Moving-patch rules}
\label{app:exp_moving}
{
The 2D advected Allen--Cahn coarsening experiment uses
the zeroth-order monitor in \eqref{eq:truncation_monitor},
with $M=0$, $u_\infty=u_+$, and $U=1$. The droplet
phase is selected by the negative-deviation branch,
\begin{equation}
    \mathcal D(t)=\{\x:u_+-u_h(\x,t)>0.45\}.
    \label{eq:patch_phase_detection}
\end{equation}
Every $100$ time steps, $u_h$ is reconstructed on a
$161\times161$ detection grid. Four-connected components
with physical area at least $0.005$ provide candidate regions.
No higher derivatives of the monitor are evaluated.

The bounding boxes are buffered and merged as described
in Section~\ref{subsec:patch_management}, with buffer
$b=0.55$ and merging gap $0.1$. The square patches in
\eqref{eq:patch_box_geometry} have minimum half-width
$h_{\min}=0.7$; centers and sizes are clipped to maintain
a distance of at least $0.2$ from the physical wall.
The reference capacity $K$ remains fixed during rollout.
The threshold $0.45$ lies between the background amplitude $0.3$ and the
droplet depth $\approx0.55$, so all droplets and no background are detected.
With decay length $\ell\approx0.12$, $b=0.55$ reduces the tail to
$\approx5\times10^{-3}$ at the patch boundary.

A patch whose detected component disappears is retained
until $\max_q|v_q|\le0.02$. If it is too close to a
proposed patch, their boxes are combined. The configuration
is rebuilt when the number of patches changes or a matched
patch's center or half-width changes by more than
$0.08h_p$, using the maximum norm for center displacement.
The old total field and background are retained during
the transfer: each new patch projects their difference
at its new nodes using \eqref{eq:patch_rebuild_local},
and the remaining field is projected onto the global
baseplate using \eqref{eq:patch_rebuild_global}.}

\section{PINN baseline}
\label{app:pinn}

The PINN (D) follows the PINN baseline of LegONet \cite{ZhangLegONet2026},
with three standard remedies so that it is not limited by avoidable training
failures. The same configuration is used in every experiment.

\paragraph{Network and constraints}
A tanh MLP $N_\theta$ with six hidden layers of width $256$ ($3.3\times10^5$
parameters) takes the normalized inputs $(\x/L,\,2t/T-1)$. The initial and
boundary conditions are imposed exactly by the output transform
\cite{Lagaris1998,LuHardPINN2021}
\begin{equation}
    v_\theta(\x,t)=v_0(\x)+\frac{t}{T}\,B(\x)\,N_\theta(\x,t),\qquad
    B(\x)=\prod_{j=1}^{d}\Big(1-\frac{x_j^2}{L^2}\Big),
\end{equation}
so that the loss contains only the PDE residual
$r_\theta=\partial_tv_\theta-\mathcal F(v_\theta)$, where $\mathcal F$ is the
right-hand side of \eqref{eq:exp_ac1d}, \eqref{eq:exp_burgers1d},
\eqref{eq:exp_acnd} or \eqref{eq:exp_acdrive} written for $v$, including the
source $s$ of the 1D Burgers forced shock and the advection of the 2D
advected Allen--Cahn coarsening.

\paragraph{Loss}
$[0,T]$ is divided into $M=32$ slabs and the slab residuals
$\mathcal L_i=\langle r_\theta^2\rangle_i$ are combined with the causal weights
\cite{WangCausalPINN2024}
\begin{equation}
    \mathcal L=\frac1M\sum_{i=1}^{M}w_i\mathcal L_i,\qquad
    w_i=\exp\Big(-\epsilon\sum_{k<i}\mathcal L_k\Big),\qquad \epsilon=1,
\end{equation}
with the weights held fixed in the gradient. Every step draws $8192$, $16384$, or $32768$ collocation points in 1D, 2D, or 3D, the same number in every slab ($512$ per slab in 2D),
from a fixed pool of $2^{20}$ (1D) or $2^{21}$ points uniform in space and
time; $v_0$ and its derivatives are precomputed on the pool.

\paragraph{Optimization and evaluation}
Adam runs for $8000$ steps with the learning rate decaying exponentially from
$10^{-3}$ to $10^{-5}$. The trained network is evaluated on the same grid and
at the same times as (A)--(C), and its error is \eqref{eq:exp_error}.

\end{document}